\documentclass[final, 3p, times]{elsarticle}

\usepackage{amssymb}
\usepackage{amsmath}
\usepackage{amsthm}

\usepackage{graphicx,color,array}
\usepackage{comment}
\usepackage{lipsum,subeqnarray}
\usepackage{amsfonts}
\usepackage{graphicx}
\usepackage{epstopdf}
\usepackage{multirow} 
\usepackage{float}
\usepackage{latexsym}
\usepackage{amsmath,bm}
\usepackage{graphicx}
\usepackage{tikz}
\usepackage{amsthm}
\usetikzlibrary{positioning,arrows.meta}
\usepackage{subfig}
\usepackage{comment}
\usepackage[utf8]{inputenc}
\usepackage{ragged2e}
\usepackage{amsmath,amssymb}
\usepackage{algorithm}
\usepackage{algorithmic}
\usepackage{placeins}
\usepackage{array}
\usepackage{adjustbox}
\usepackage[utf8]{inputenc}
\usepackage[numbers]{natbib}

\usepackage{subeqnarray}
\usepackage{graphicx}

\usepackage{lineno}
\usepackage{hyperref}
\usepackage{xcolor}
\usepackage{booktabs}

\journal{Journal of Computational Physics}

\begin{document}

\begin{frontmatter}



\title{Green Multigrid Network}  

\author[JLU]{Ye Lin}
\ead{linye21@mails.jlu.edu.cn}
\affiliation[JLU]{
organization={School of Mathematics, Jilin University}, 
addressline={2699 Qianjin Street}, 
city={Changchun},
postcode={130012}, 
state={Jilin},
country={China}.}

\author[TXST]{Young Ju Lee}
\ead{yjlee@txstate.edu}
\affiliation[TXST]{
organization={Department of Mathematics, Texas State University}, 
addressline={601 University Drive}, 
city={San Marcos},
postcode={78666}, 
state={Texas},
country={U.S.A.}.}

\author[JLU]{Jiwei Jia\corref{*}}
\ead{jiajiwei@jlu.edu.cn}
\cortext[*]{Corresponding Author.}


\begin{abstract}

GreenLearning networks (GL)\cite{boulle2022gl} directly learn Green's function in physical space, making them an interpretable model for capturing unknown solution operators of partial differential equations (PDEs). For many PDEs, the corresponding Green's function exhibits asymptotic smoothness. In this paper, we propose a framework named Green Multigrid networks (GreenMGNet), an operator learning algorithm designed for a class of asymptotically smooth Green's functions.

Compared with the pioneering GL, the new framework presents itself with better accuracy and efficiency, thereby achieving a significant improvement. GreenMGNet is composed of two technical novelties. First, Green's function is modeled as a piecewise function to take into account its singular behavior in some parts of the hyperplane. Such piecewise function is then approximated by a neural network with augmented output(AugNN) so that it can capture singularity accurately. Second, the asymptotic smoothness property of Green's function is used to leverage the Multi-Level Multi-Integration (MLMI) algorithm for both the training and inference stages. Several test cases of operator learning are presented to demonstrate the accuracy and effectiveness of the proposed method. On average, GreenMGNet achieves $3.8\%$ to $39.15\%$ accuracy improvement. To match the accuracy level of GL, GreenMGNet requires only about $10\%$ of the full grid data, resulting in a $55.9\%$ and $92.5\%$ reduction in training time and GPU memory cost for one-dimensional test problems, and a $37.7\%$ and $62.5\%$
reduction for two-dimensional test problems. 
\end{abstract}



\begin{keyword}
 Neural network \sep Operator learning \sep Green's function \sep Multigrid Method
\end{keyword}

\end{frontmatter}


\section{Introduction} \label{sec:intro}

Over the past few years, deep neural networks have been recognized as a potent technique for solving partial differential equations (PDEs) across a broad range of scientific and engineering problems \cite{karniadakis2021physics}. The approaches within this field can be roughly divided into two main categories: (1) single PDE solvers and (2) operator learning. Single PDE solvers, such as physics-informed neural networks (PINNs) \cite{raissi2019physics, cuomo2022scientific, wang2023expert}, the deep Galerkin method \cite{sirignano2018dgm}, the deep Ritz method \cite{yu2018deep}, and finite neuron method (FNM) \cite{xu2020fnm}, optimize deep neural network by the PDE residual or variational energy formulation of the PDE to obtain solution function. These methods are tailored to solve a particular instance of the PDE. Thus, if there are changes in the initial condition, boundary condition, or forcing term, it is unavoidable to retrain the neural network and it is computationally inefficient in practice. 

On the other hand, Operator Learning employs neural networks to learn operators between function spaces, which could be utilized to learn solution operators of PDEs from data pairs. Several operator learning methods have been recently proposed, such as deep green networks (DGN) \cite{boulle2022gl, gin2021deepgreen, lin2023bi}, deep operator network (DeepONet) \cite{lu2021learning, wang2021learning, jin2022mionet}, neural operators (NOs) \cite{li2020gno, li2020fourier, rahman2023uno, tripura2023wavelet, he2024mgno}. Despite the inspirations behind these methods are distinct, they all establish a connection between operator learning and integral transform in different forms, which are closely related to Green's function methods. 

For linear PDEs, the solution operator can be represented as an integral transform with a Green's function, which is typically unknown. It is natural to model this unknown function with a neural network and learn it from data. In \cite{boulle2022gl}, this idea is pioneered to design GreenLearning networks (GL) for solving operator learning tasks. However, there is room for improvement, both in approximation accuracy and efficiency for GL. First, Green's functions, in general, possess singularities, such as those for elliptic PDEs \cite{borm2003introduction, bebendorf2003existence, boulle2023ellipticrandSVD}. Thus, a direct neural network approximation of Green's function with singularities can present itself a difficulty in obtaining a good approximation property \cite{boulle2020rational, hu2022discontinuity}. Second, the numerical approximation of integral transform can be written as a dense matrix-vector product\cite{boulle2023review}. Such a dense matrix-vector product can be a computational bottleneck in GL. For instance, we assume that the domain $\Omega$ of the target problem is discretized on an equidistant grid with a total of $n$ nodes, then the evaluation of the dense kernel matrix, the discrete Green's function, requires $O(n^2)$ neural network forward operations, and the calculation of the dense matrix-vector product itself will also require $O(n^2)$ operations.

There are several works on improving neural network approximation ability for discontinuous function fitting, such as those with sharp edges, cusps, and singularities. Since discontinuities are common in the solutions of interface problems, many contributions originate from related research fields. \cite{hu2022discontinuity,tseng2023cusp} add the level set function of the interface as an additional input to the network to retain the cusps in the solution. \cite{wu2022inn, wu2023solving} decomposes the computational domain into several subdomains and fits solutions by different networks to keep discontinuities across interfaces. Nevertheless, these techniques can not be directly used to learn kernel functions. On the other hand, numerous effective techniques have been devised to reduce the complexity of dense matrix-vector product calculation, such as Fast Fourier Transform (FFT) \cite{cooley1965algorithm} and fast multi-resolution algorithms \cite{greengard1987fast,brandt1990multilevel,beylkin1991fast}. Several fast methods have already been applied in operator learning algorithms \cite{li2020gno, li2020multipole, gupta2021multiwaveletbased}. However, their application to GL setting is missing in the literature.

This paper is aiming at filling these gaps. The main contribution of the proposed work is the development of Green Multigrid networks (GreenMGNet). The main technical novelties in our work are motivated by previous works on interface problems and fast matrix-vector product algorithms. 
GreenMGNet is designed to effectively learn kernels under asymptotically smoothness assumption. More precisely, we overcome the approximation difficulties of the neural network for Green's function due to the singularity via a splitting technique of the neural network. Namely, we split the domain of Green's function into two parts around the hyperplane $\sum_{i=1}^d(x_i - y_i) = 0$, where $\mathbf{x} = (x_1, x_2, \cdots x_d)$, $\mathbf{y} = (y_1, y_2, \cdots y_d)$, and $d$ is the dimension. We then model the kernel function as a piece-wise smooth function as done for the interface problem. We note that splitting domain around $\sum_{i=1}^d(x_i - y_i) = 0$ will introduce extra non-smoothness. This is remedied by applying an aggregation strategy to solve this problem. Finally, an augmented output neural network (AugNN) is designed to fit the piecewise kernel function, which avoids the explicit discontinuity in neural network training. Secondly, under the asymptotic smoothness assumption of the underlying kernels, the corresponding matrix-vector product can be computed efficiently by employing the Multi-Level Multi-Integration (MLMI)\cite{brandt1990multilevel} algorithm. Note that such an algorithm is originally designed to solve integral equations with potential-type kernels. Utilizing MLMI during the GL training stage allows us to estimate the dense matrix-vector product using a subset of full grid points, which significantly reduces the operations of neural networks for kernel evaluation. During the inference stage, no neural network forward operation is necessary, as the discrete kernel can be computed in advance. Consequently, the computational cost arises solely from discrete kernel integration, where MLMI proves more efficient than direct matrix-vector multiplication.

Our main contributions can be summarized as follows: 
\begin{itemize}
\item[(1)] We model Green's function estimation problem as an interface problem, and design AugNN to reduce the learning difficulties for Green's function with diagonal discontinuity. 
\item[(2)] We design GreenMGNet by leveraging the MLMI and AugNN, which evaluates discrete kernel integrals with fewer points, reducing data points requirement without worsening the performance.
\item[(3)] We test GreenMGNet by learning Green's function from various elliptic PDEs. It is shown that GreenMGNet surpasses GL in terms of both accuracy and efficiency for all the test problems.
\end{itemize}
Even if the GreenMGNet is mainly applied to obtain the Green function for the elliptic PDEs, the framework can be equally well applied for learning any kernel functions with asymptotically smoothness property through integral transform. 

The rest of the paper is organized as follows. In Section \ref{sec:ol_n_gl}, we overview the basics of the operator learning problem and GL method. In Section \ref{sec:ask_n_augnn}, we introduce the asymptotically smooth kernel and the design of AugNN for kernel fitting in such cases. In section \ref{sec:mlmi_n_gmgn}, we present the MLMI algorithm and GreenMGNet method for operator learning. Several numerical experiments are used to verify the effectiveness of both AugNN and GreenMGNet in section \ref{sec:num_exps}, and finally a conclusion is given in section \ref{sec:conclusion}.

\section{Operator Learning and GreenLearning} \label{sec:ol_n_gl}

\subsection{Operator Learning} \label{subsec:ol}
Let $\mathcal{U}$ and $\mathcal{V}$ be Banach spaces of functions defined on $d$-dimensional spatial domain $\Omega \subset \mathbb{R}^d$, and $\mathcal{A} : \mathcal{U} \rightarrow \mathcal{V}$ be an operator which maps between $\mathcal{U}$ and $\mathcal{V}$. Suppose we have observations $\{f^{(k)}, u^{(k)}\}^{N}_{k=1}$, where $f^{(k)} \in \mathcal{U}$ and $u^{(k)} \in \mathcal{V}$, the aim of operator learning is to learn a model $\mathcal{A}_{\theta}$ with a finite parameters $\theta \in \mathbb{R}^p$ to approximate $\mathcal{A}$ by minimizing:
\begin{equation}
    \label{eq:operator_learning}
    \min_{\theta} \frac{1}{N} \sum_{k=1}^{N} C(\mathcal{A}_{\theta}(f^{(k)}), u^{(k)}),
\end{equation}
where $C: \mathcal{V} \times \mathcal{V} \rightarrow \mathbb{R}$ is a cost functional, such as relative $L^2$ norm, which measures the discrepancy between $\mathcal{A}_{\theta}(f)$ and $u$. Typically, we select a neural network as $\mathcal{A}_{\theta}$, where $\theta$ represents the weights and biases of the neural network.

\subsection{GreenLearning}
\begin{figure}[h]
\centering
\includegraphics[scale=.7]{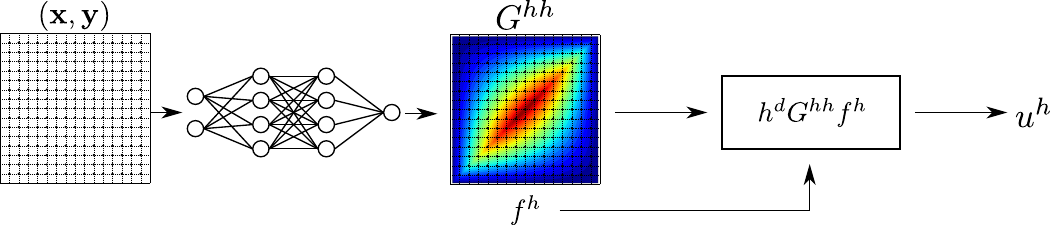}
\caption{Schematic of Green Learning}
\label{fig:gl-arch}
\end{figure}

GreenLearning networks (GL)\cite{boulle2022gl} assume the operator $\mathcal{A}$ can be represented as an integral operator with a kernel $G:\Omega \times \Omega \rightarrow \mathbb{R} \cup \{\infty\}$ as \cite{evans2022partial, boulle2022gl, boulle2023review, boulle2023ellipticefficient}
\begin{equation} 
\label{eq:green_func}
\mathcal{A}(f)(\mathbf{x}) = u(\mathbf{x}) = \int_{\Omega} G(\mathbf{x},\mathbf{y}) f(\mathbf{y}) \mathrm{d} \mathbf{y}, \quad \mathbf{x} \in \Omega,
\end{equation}
and learn $G$ by a neural network $G_{\theta}$ from data pairs, where $\theta$ denotes the learnable parameters of a neural network. Numerical evaluation of (\ref{eq:green_func}) is necessary, so the GL approach can be divided into two steps as shown in Figure \ref{fig:gl-arch}, the first is discrete kernel function evaluation by a neural network, the second is kernel integral estimation by some quadrature rule.

To be specific, assume $\mathbf{x}_i=\mathbf{x}_0 + ih$ are equidistant grid points defined on a $d$-dimensional domain $\Omega$, where $i=(i_1, i_2, \cdots, {i_d} )$ is a vector of integers and $h$ is the mesh size. In the same way, we define grid points $\mathbf{y}_j$. Let $u^h$ and $f^h$ are the discrete functions approximating $u$ and $f$ on the grids, and we denote $u^h_i = u(\mathbf{x}_i)$ and $f^h_j=f(\mathbf{y}_j)$, respectively. Using simple midpoint quadrature rule, the integral (\ref{eq:green_func}) are discretized as:
\begin{equation} 
\label{eq:multi-sum}
u^h_i = h^d \sum_{j} G^{hh}_{i,j} f^h_j,
\end{equation}
where $G^{hh}_{i,j} = G_{\theta}(\mathbf{x}_i, \mathbf{y}_j)$ is a neural network approximation of $G$ on $(\mathbf{x}_i, \mathbf{y}_j)$. We denote $G^{hh}$ as the discrete kernel matrix evaluated on the full grid point pairs, then the matrix form of (\ref{eq:multi-sum}) becomes $u^h = h^d G^{hh} f^h$. Given dataset $\{ (f^h)^{(k)}, (u^h)^{(k)} \}^{N}_{k=1}$, neural network $G_{\theta}$ can be trained by minimizing the  relative $L^2$ loss:
\begin{equation} 
\label{eq:rl2}
\text{Loss} = \frac{1}{N} \sum_{k=1}^{N} \frac{\| (u^h)^{(k)} - h^d (G^{hh}) (f^h)^{(k)} \|}{ \| 
(u^h)^{(k)} \| }  
\end{equation}

GL directly learns the kernel(Green's function) in the physical space. This contrasts with other operator learning methods like DeepONets and FNOs, which are designed to learn the actions of operators. Some approaches \cite{lin2023bi,peng2023deep,sun2023binn,negi2024learning,mod-net} also recover Green's functions using a neural network, but they need a PDE loss for optimization, similar to PINNs.

Several obvious issues with the GL method need to be addressed in practice. First, the GL method assumes that the solution operator of PDEs can be written as a kernel integral. However, this assumption is not appropriate for nonlinear PDEs. Second, approximating functions with cusps or singularities, which are common for Green's function, may be difficult for a standard neural network. Additionally, using a typical quadrature rule to evaluate kernel integrals results in high computational costs. 

To improve the accuracy, \cite{gin2021deepgreen, mod-net} introduced an additional neural network as a nonlinear transformation after the kernel integral calculation. To reduce learning difficulties, \cite{boulle2022gl} employed the Rational neural network\cite{boulle2020rational} due to its powerful approximation capabilities, particularly for approximating unbounded or singular Green's function. Even so, a Rational neural network usually gives a smooth approximation where cusps should exist. Another approach from \cite{lin2023bi} decomposes Green's function into a singular part and a smooth part when the exact location of singularities is known. However, determining the exact location of singularities remains challenging. For the efficiency issue, to the best of our knowledge, no fast algorithm has been implemented to improve the efficiency of GL. Nevertheless, \cite{boulle2023ellipticefficient, boulle2023ellipticrandSVD, boulle2022parabolic} provide some theoretical results on the amount of training data pairs requirement for elliptic and parabolic problems, which indirectly reduce training cost on GL.

In this paper, we focus on solving the last two problems mentioned above. We propose an augmented output network to further improve the approximation accuracy and apply a fast matrix-vector product algorithm to reduce the computational cost of GL.

\section{Asymptotically smooth kernel and augmented output neural network} \label{sec:ask_n_augnn}

\subsection{Asymptotically smooth kernel}

In many applications, kernels $G(\mathbf{x}, \mathbf{y})$ are not smooth over the entire domain. They exhibit unbounded derivatives as $\mathbf{y} \rightarrow \mathbf{x}$. On the other hand, the smoothness of $G$ increases indefinitely with the distance $\| \mathbf{y} - \mathbf{x} \|$. Generally, we refer to $G(\mathbf{x}, \mathbf{y})$ asymptotically smooth with respect to both $\mathbf{y}$ and $\mathbf{x}$ if, for any positive integer $q$ and $q$-order derivative $\partial^{q}$, 
\begin{equation}
\label{eq:asmooth-y} 
    \ \| \partial^{q} G(\mathbf{x},\mathbf{y}) \| \leq  C_{q} {\|\mathbf{x}-\mathbf{y}\|}^{g-q},    
\end{equation}
where $g$ is independent of $q$, $C_{q}$ depends only on $p$\cite{brandt1991multilevel}. An example is potential-type kernels, such as Green's function for the two-dimensional infinite-domain Poisson problem $G(\mathbf{x}, \mathbf{y}) = \mathrm{ln}(\|\mathbf{x}-\mathbf{y}\|)/(2\pi)$ and three-dimensional infinite-domain Poisson problem $G(\mathbf{x}, \mathbf{y}) = -1 / (4\pi \|\mathbf{x}-\mathbf{y}\|)$. One should note that, for uniformly elliptic PDEs, the corresponding Green's function is not as general as asymptotically smooth kernel, typically they also have the low-rank property and a slower decay limitation\cite{boulle2023ellipticrandSVD, bebendorf2003existence, gruter1982green}.


For asymptotically smooth kernels, singularities may only occur at the diagonal and the function is smooth everywhere else\cite{borm2003introduction, brandt1990multilevel, brandt1991multilevel}. Singularity is a type of non-smoothness that can not be easily handled by a regular neural network. Most neural networks are inherently smooth due to the corresponding smooth activation function. In most cases, we can numerically approximate a singularity with a bounded cusp which is also not smooth. Inspired by techniques recently used in elliptic interface problems \cite{tseng2023cusp,hu2022discontinuity, wu2022inn,he2022mesh}, we treat an asymptotically smooth kernel as a piecewise function, which could capture the cusp more accurately. Many of Green's functions are asymptotically smooth, such as uniformly elliptic PDEs. To illustrate the ideas, we discuss Green's function of the Poisson equation.

\begin{figure}[h]
\centering
\includegraphics[scale=.7]{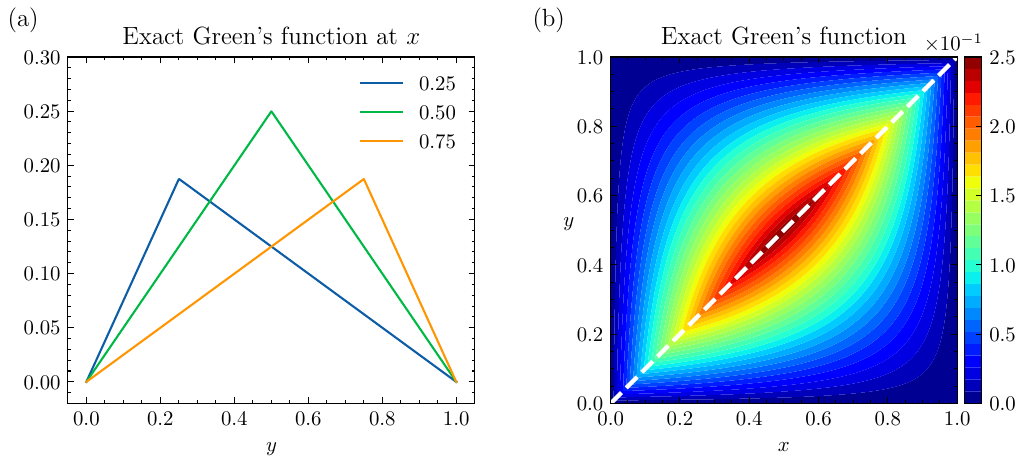}
\caption{(a): The slice of Green's function at $x=0.25,x=0.5,x=0.75$ and the cusp appears at $x=y$. (b): The Green's function on $[0,1] \times [0,1]$ and $x=y$(white dash line) divides the domain into two subdomains as an interface.}
\label{fig:poisson1d_green}
\end{figure}

\paragraph{One-dimensional Poisson equation}\label{para:poisson_1d} We first give the example of one-dimensional Poisson equation with homogeneous boundary conditions:

\begin{equation} 
\label{eq:poisson1d}
\begin{split}
-\frac{d^2 u}{d x^2} & = f(x), \quad x \in [0,1] \\
u(0) & = u(1) = 0.
\end{split}
\end{equation}
The associated Green's function is available in closed-form\cite{li2020gno} as 
\begin{equation} 
\label{eq:poisson1d_green}
\begin{split}
G(x,y) & = \frac{1}{2}(x+y-|y-x|) - xy, \quad x,y \in [0,1], \\
\end{split}
\end{equation}
or we can rewrite (\ref{eq:poisson1d_green}) as a piecewise function:
\begin{equation} 
\label{eq:poisson1d_green_piecewise}
G(x,y) = \begin{cases} 
      x(1-y), & x-y\leq 0, \\
      y(1-x), & x-y>0.  
   \end{cases}
\end{equation}
In (\ref{eq:poisson1d_green_piecewise}), $x-y=0$ acts as an interface splitting the domain into two subdomains, as depicted in Figure \ref{fig:poisson1d_green}. Although singularity doesn't exist in this case, (\ref{eq:poisson1d_green_piecewise}) is still an asymptotically smooth kernel, and the presence of the cusp at $x=y$ still increases the difficulty for a neural network approximation. 

\begin{figure}[h]
\centering
\includegraphics[scale=.7]{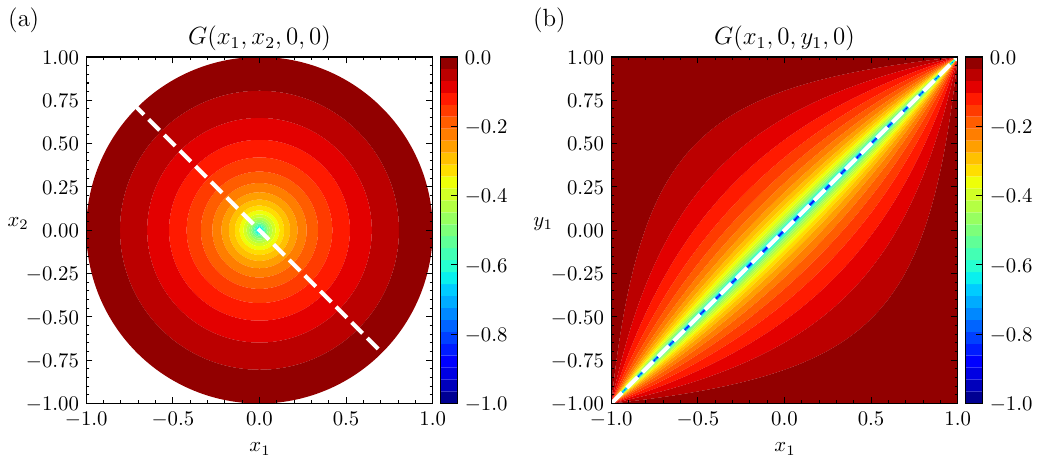}
\caption{(a): the slice of Green's function at $(x_1, x_2, 0, 0)$, where singularity exists at $(0,0,0,0)$. (b): the slice of Green's function at $(x_1, 0, y_1, 0)$ where singularity exists at $x_1=y_1$. The white dash line is the hyperplane $x_1 - y_1 + x_2 - y_2 = 0$ on each slice.}
\label{fig:poisson2d_green}
\end{figure}

\paragraph{Two-dimensional Poisson equation} \label{para:poisson_2d} For two-dimensional Poisson equation with homogeneous boundary condition on the unit disk:
\begin{equation} 
\label{eq:poisson2d}
\begin{split}
-\nabla \cdot (\nabla u(\mathbf{x})) &= f(\mathbf{x}), \quad \mathbf{x} \in \Omega, \\
u(\mathbf{x}) &= 0, \quad \mathbf{x} \in \partial \Omega,
\end{split}
\end{equation}
where $\Omega=D(0,1)$.
In this case, the closed form of the associated Green's function is known and given by \cite{myint2007linear}:
\begin{equation} 
\label{eq:green-disk-poisson} 
G(\mathbf{x}, \mathbf{y}) =  \frac{1}{4 \pi} \mathrm{ln} \left ( \frac{(x_1-y_1)^2 + (x_2-y_2)^2}{(x_1 y_2 - x_2 y_1)^2 + (x_1y_1 + x_2y_2 - 1)^2} \right ), 
\end{equation} 
where $\mathbf{x} = (x_1,x_2), \mathbf{y} = (y_1,y_2)$ are two points in the unit disk. The Green's function (\ref{eq:green-disk-poisson}) is an asymptotically smooth kernel with singularity at $\mathbf{x} = \mathbf{y}$. Unlike the one-dimensional case, (\ref{eq:green-disk-poisson}) is a function defined in a four-dimensional domain that can not be separated by $\mathbf{x}=\mathbf{y}$. We apply a hyperplane $x_1 - y_1 + x_2 - y_2 = 0$ to split the four-dimensional domain into two parts(See Figure \ref{fig:poisson2d_green}) and rewrite Green's function in a piecewise form as follows:
\begin{equation} 
\label{eq:green-2d-piecewise} 
    G(\mathbf{x}, \mathbf{y}) =  \begin{cases} 
      G_1(\mathbf{x}, \mathbf{y}), & x_1-y_1+x_2-y_2 \leq 0, \\
      G_2(\mathbf{x}, \mathbf{y}), & x_1-y_1+x_2-y_2 > 0.  
   \end{cases} 
\end{equation} 

Generally, for a $d$-dimensional uniformly elliptic problem, the corresponding Green's function $G(\mathbf{x}, \mathbf{y})$ is a $2d$-dimensional function, and singularity may only exist at $\mathbf{x}=\mathbf{y}$. Instead of directly fitting $G(\mathbf{x}, \mathbf{y})$, we consider it as a piecewise function:
\begin{equation} 
\label{eq:green-nd-piecewise} 
    G(\mathbf{x}, \mathbf{y}) =  \begin{cases} 
      G_1(\mathbf{x}, \mathbf{y}), & \sum_{i=1}^d(x_i - y_i) \leq 0, \\
      G_2(\mathbf{x}, \mathbf{y}), & \sum_{i=1}^d(x_i - y_i) > 0.  
   \end{cases} 
\end{equation} 
where $\sum_{i=1}^d(x_i - y_i) = 0$ is the splitting interface. The piecewise function (\ref{eq:green-nd-piecewise}) can be fitted separately in each subdomain to preserve non-smoothness on the interface. 

\subsection{Augmented output neural network}\label{sec:aug-nn}
GL method uses a neural network with a single output to model a scalar output kernel. A typical neural network is good enough for approximation when the underlying kernel is smooth. For cases we discussed above, the approximation error will be large around the cusp or singularity because of non-smoothness(See numerical results in Section \ref{subsec:aug_nn-result}). To solve this problem, we proposed an augmented output neural network(AugNN) to model piecewise function, which gives a better approximation near non-smoothness. 

\begin{figure}[h]
\centering
\includegraphics[scale=.75]{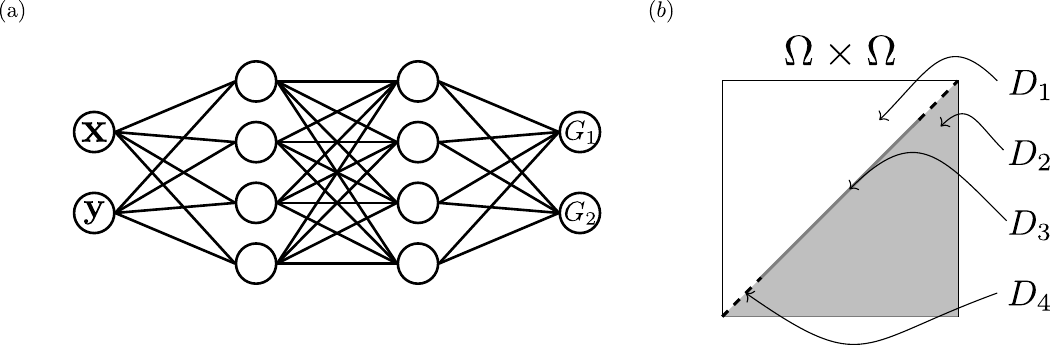}
\caption{(a): the architecture of AugNN. (b): Illustration of subdomains introduced in (\ref{eq:subdomains}), averaging $G_1$ and $G_2$ is need on $D_4$ for smoothing.}
\label{fig:augnn-domain}
\end{figure}

The structure of AugNN is shown in Figure \ref{fig:augnn-domain}. For a scalar output kernel, AugNN has outputs $G_1$ and $G_2$ as defined in (\ref{eq:green-nd-piecewise}). Despite non-smoothness is kept on the interface, extra non-smoothness is introduced since $\mathbf{x} = \mathbf{y}$ is a subset of $\sum_{i=1}^d(x_i - y_i) = 0$ when $d > 1$. For instance, the slice of Green's function $G(x_1,x_2,0,0)$ shown in Figure \ref{fig:poisson2d_green} possesses a singularity only at the point $(0,0,0,0)$. On this slice, the interface, which is assumed to have non-smoothness, forms a line. However, excluding kernel value at $(0,0,0,0)$, Green's function should be smooth along this line. 

To solve this problem, the domain of $G$ is divided into four subdomains, defined as:
\begin{equation}
\label{eq:subdomains}
\begin{split}
    D_1 &= \{(\mathbf{x},\mathbf{y}) | \sum_{i=1}^d(x_i - y_i) < 0\},\\
    D_2 &= \{(\mathbf{x},\mathbf{y}) | \sum_{i=1}^d(x_i - y_i) > 0\},\\
    D_3 &= \{(\mathbf{x},\mathbf{y}) | \mathbf{x} = \mathbf{y} \},\\            
    D_4 &= \{(\mathbf{x},\mathbf{y}) | (\sum_{i=1}^d(x_i - y_i) = 0) \backslash (\mathbf{x} = \mathbf{y}) \},\\
\end{split}
\end{equation} 
$D_1$ and $D_2$ are points on the two sides of the interface, and $D_3$ and $D_4$ are points on the interface. Specifically, $D_3$ lies on the diagonal where non-smoothness may exist and $D_4$ are 
points also on the interface where the values of kernel function are expected to be smooth. We apply a simple aggregation strategy to maintain the smoothness on $D_4$ by averaging $G_1$ and $G_2$:
\begin{equation} 
\label{eq:augnn-piecewise} 
    G(\mathbf{x}, \mathbf{y}) =  \begin{cases} 
      G_1(\mathbf{x}, \mathbf{y}), & (\mathbf{x}, \mathbf{y}) \in D_1, \\
      G_2(\mathbf{x}, \mathbf{y}), & (\mathbf{x}, \mathbf{y}) \in D_2, \\
      G_2(\mathbf{x}, \mathbf{y}), & (\mathbf{x}, \mathbf{y}) \in D_3, \\
      \frac{1}{2}(G_1(\mathbf{x}, \mathbf{y})+G_2(\mathbf{x}, \mathbf{y})), & (\mathbf{x}, \mathbf{y}) \in D_4, \\
   \end{cases} 
\end{equation} 
We take the average of $G_1$ and $G_2$ to remove the undesired non-smoothness. By doing this, the augmented output neural network keeps non-smoothness only at $\mathbf{x}=\mathbf{y}$, which makes them a better model for asymptotically smooth kernel learning. We remark that instead of taking the average, one can introduce some penalty loss such as 
\begin{equation}
\label{eq:smooth_loss} \sum_{D_4} \| G_1(\mathbf{x}, \mathbf{y}) - G_2(\mathbf{x}, \mathbf{y}) \|,
\end{equation}
which plays a role in achieving continuity in a weak sense as in the Discountinuous Galerkin method \cite{arnold1982interior}. However, we take a simple path to take the average in this domain $D_4$, which seems to present no issue in accuracy. 

\section{Multi-Level 
Multi-Integration and GreenMGNet} \label{sec:mlmi_n_gmgn}

\subsection{Multi-Level Multi-Integration}
Discretizing the integral transform (\ref{eq:green_func}) using a quadrature rule, it will become a matrix-vector product as (\ref{eq:multi-sum}). If the matrix has arbitrary entries, no way can exist which is significantly faster than the straightforward multiplication. In most applications of interest, however, the matrix has a certain structure, which can be circulant, low-rank, banded, or hierarchical off-diagonal low-rank(HODLR) and reflect the characteristics of
the corresponding kernel function \cite{boulle2023review}. 

Based on the structure prior, many fast matrix-vector product algorithms have been developed, such as Fast Fourier Transform(FFT) \cite{cooley1965algorithm}, Fast Wavelet Transform(FWT) \cite{beylkin1991fast}, Fast Multipole Method(FMM) \cite{greengard1987fast} and Multi-level Multi-Integration(MLMI)\cite{brandt1990multilevel}. Many fast methods have already been used in operator learning tasks by assuming the underlying properties of kernels. Fourier neural operator(FNO) \cite{li2020fourier} and its variants \cite{tran2023factorized,li2023fourier,you2022learning} employ FFT in the design of the network architecture. DeepONet \cite{lu2021learning} and its variants \cite{wang2021learning,jin2022mionet,wang2022improved} are motivated by universal operator approximation theorem\cite{chen1995universal}, which is a form of low-rank approximation of operators. Another architecture related to low-rank approximation is the attention block\cite{vaswani2017attention}, which has been utilized in the design of many operator learning algorithms\cite{cao2021choose,hao2023gnot,kissas2022learning,guo2022transformer,lee2022meshindependent}. Multipole graph neural operator(MGNO)\cite{li2020multipole} is designed with the inspiration from Fast Multipole Methods(FMM). \cite{gupta2021multiwaveletbased} introduces a multiwavelet-based operator by parametrizing the integral operators using Fast Wavelet Transform(FWT). However, these methods learn the action of the operator rather than kernel functions like GL. As far as we know, no research has yet incorporated fast matrix-vector product algorithms to enhance the efficiency of GL. 

In this paper, our interest is in learning unknown asymptotically smooth kernels, where MLMI can be used for fast matrix-vector product calculation. MLMI method was proposed in \cite{brandt1990multilevel} for fast evaluation of the integral transforms involving both smooth and asymptotically smooth kernels. Some extensions and applications of MLMI can be found in \cite{brandt1991multilevel, grigoriev2004fast, polonsky1999numerical, polonsky2000fast, venner2000multi}. In this section, we first introduce the algorithm of MLMI, then the Green Multigrid networks(GreenMGNet) is presented which leverages MLMI to improve the efficiency of GL.

\paragraph{Notation of MLMI}
MLMI is designed to obtain the numerical evaluation of the integral transforms like (\ref{eq:green_func}). The full algorithm uses a target grid with mesh size $h$ and a series of coarser grids recursively. For simplicity, we assume $\Omega = [0,1]$, and only two grids are involved. A target fine grid is defined as $x^h_i = ih$, where $i \in [1,2,\cdots,n]$ is the index and $h=1/(n+1)$ is the meshsize. The value of function $u$ at a fine grid $i$ with location $x_i$ denotes $u^h_i$. Similarly, an auxiliary coarse grid is defined as $X^H_I = IH$, where $I \in [1,2,\cdots,N]$ and $H=1/(N+1)$, and $U^H_I$ is the value of $U$ at a coarse grid point $I$ with location $X_I$. We assume $H=2h$ for simplicity, so $n=2N+1$ and satisfying $x^H_I = x^h_{2J}$. For discrete kernel, $G^{hh}_{i,j}$ approximates $G(x^h_i, x^h_j)$, and its index may from different grid, such as $G^{hH}_{i,J}$ approximates $G(x^h_i, x^H_J)$.

We will use $I^{h}_{H}$ to denote an interpolation operator from coarse grid to fine grid; i.e $I^h_H v^H$ stands for a coarse grid$(H)$ function $v^H$ is interpolated to a fine grid$(h)$. For each fine grid index $i$, $I^h_{H} G^{hH}_{i,\cdot}$ denotes a fine-grid function obtained by interpolating from the coarse-grid function $G^{hH}_{i,\cdot}$, and the value of the interpolated function at $x^h_j$ is denoted by $[I^h_{H} G^{hH}_{i,\cdot}]_j$. The adjoint of an interpolation operator is a restriction operator $R^{H}_{h}$, such as $R^H_h v^h$ represents a coarse grid function which is restricted from a fine grid function $v^h$, and $[R^H_h v^h_{\cdot}]_J$ is the value of the restricted function at $x^H_J$. We only consider linear interpolation if $d=1$ or bilinear if $d=2$ in this paper. If we deem $v^H$ as a $N \times 1$ vector, the interpolation operator $I^h_{H}$ can be written as an $n \times N$ matrix, then restriction operator is given as $R^H_{h} = 2^{-d}(I^h_H)^T$.

In practice, we will use more than two grids. It is convenient to refer to them as levels, starting with the coarsest grid that will be level 1, the next finer grid being level 2, etc.

  


\paragraph{Smooth kernels}
Whenever the kernel $G(x,y)$ is sufficiently smooth with respect to $y$, we can consider the following interpolation:
\begin{equation} 
\label{eq:interp}
\Tilde{G}^{hh}_{i,j} = [I^h_H G^{hH}_{i,\cdot}]_j,\\
\end{equation}
where $I^h_H$ is the interpolation operator and $G^{hH}_{i,\cdot}$ is injected from $G^{hh}_{i,\cdot}$; i.e., $G^{hH}_{i,J} := G^{hh}_{i,2J}$. We then note that (\ref{eq:multi-sum}) can be approximated by

\begin{equation} 
\label{eq:multi-sum_sim}
\begin{split}
u^h_i &\simeq  h \sum_{j} \Tilde{G}^{hh}_{i,j} f^h_j = h \sum_{j} [I^h_H G^{hH}_{i,\cdot}]_j f^h_j \\
    &=h \sum_{j}  G^{hH}_{i,J} [R^h_H f^h_{\cdot}]_J = H \sum_{J} G^{hH}_{i,J} \Tilde{f}^H_J.
\end{split}
\end{equation}

Whenever $G(x,y)$ is also sufficiently smooth as a function of $x$, the value of $u^h_i$ can be calculated only for coarse grid points $i=2I$, using the interpolation to obtain the other values on the fine grid. Namely, we have 
\begin{equation} 
\label{eq:smooth_sim}
u^h_i \simeq I^h_H \Tilde{u}^H_I,
\end{equation}
where 
\begin{equation} 
\label{eq:TuH}
\Tilde{u}^H_I := H \sum_{J} G^{HH}_{I,J} \Tilde{f}^H_J,
\end{equation}
and where $G^{HH}_{\cdot, J}$ is "injected" from $G^{hH}_{\cdot,J}$, i.e., $G^{HH}_{I,J} := G^{hH}_{2I,J}=G^{hh}_{2I,2J}$. 

The matrix-vector product (\ref{eq:multi-sum}) on the fine grid has thus been reduced to the analogous coarse grid matrix-vector product (\ref{eq:TuH}) with interpolation (\ref{eq:smooth_sim}). The latter problem can be coarsened repeatedly until the number of gridpoints is proportional to $\sqrt{n}$. Then, the complexity is reduced from $O(n^2)$ to $O(n)$ in this case\cite{brandt1990multilevel}\cite{mlmi2005acoustics}. 


\paragraph{Asymptotically smooth kernels} If $G(x,y)$ is asymptotically smooth, it has some singular points, and the smoothness increases rapidly with increasing distance from these points. For simplicity, we also assume that the only singular points are points on diagonals, i.e. $x = y$.

Let's start with the case that $i=2I$ and derive an exact expression relating $\Tilde{u}^H_I$ and $u^h_i$, replacing the approximate equation (\ref{eq:multi-sum_sim}),
\begin{equation} 
\label{eq:even-decomp}
\begin{split}
u^h_i &= h \sum_{j} G^{hh}_{i,j} f^h_j = h \sum_{j} \Tilde{G}^{hh}_{i,j} f^h_j + h \sum_{j} (G^{hh}_{i,j} - \Tilde{G}^{hh}_{i,j}) f^h_j \\
    &=\Tilde{u}^H_I + h \sum_{j} (G^{hh}_{i,j} - \Tilde{G}^{hh}_{i,j}) f^h_j. \\
\end{split}
\end{equation}
(\ref{eq:even-decomp}) splits $u^h_i$ into a smooth term $\Tilde{u}^H_I$ and a correction term $h \sum_{j} (G^{hh}_{i,j} - \Tilde{G}^{hh}_{i,j}) f^h_j$. Using the knowledge of the differences $(G^{hh}_{i,j} - \Tilde{G}^{hh}_{i,j})$ will be small when $\|i-j\| \gg m$\cite{brandt1990multilevel}, where $m$ is some threshold represents the local correction range, the correction term can be approximated by points near the diagonal,
\begin{equation} 
\label{eq:even-cor_local}
\Tilde{e}^h_{i} := \sum_{|j-i| \leq m} (G^{hh}_{i,j} - \Tilde{G}^{hh}_{i,j}) f^h_j, \\
\end{equation}
so that
\begin{equation} 
\label{eq:even-cor_approx}
u^h_i \simeq \Bar{u}^H_I = \Tilde{u}^H_I + \Tilde{e}^h_{2I}. \\
\end{equation}

If the point $i$ is not a coarse grid point$(i \neq 2I)$, we define another coarse grid approximation,
\begin{equation} 
\label{eq:interp_odd}
\Hat{G}^{hh}_{i,j} = [I^h_H G^{Hh}_{\cdot, j}]_i,\\
\end{equation}
where $G^{Hh}_{I,j} = G^{hh}_{2I,j}$. Again, we split $u^h_i$ into a smooth term and a correction term as (\ref{eq:even-decomp}), 
\begin{equation} 
\label{eq:odd_decomp}
u^h_i = h \sum_{j} G^{hh}_{i,j} f^h_j = h \sum_{j} \Hat{G}^{hh}_{i,j} f^h_j + h \sum_{j} (G^{hh}_{i,j} - \Hat{G}^{hh}_{i,j}) f^h_j.
\end{equation}
The smooth term $ h \sum_{j} \Hat{G}^{hh}_{i,j} f^h_j$ can be approximated by interpolating $\Bar{u}^H_I$, 
\begin{equation} 
\label{eq:odd_smooth_approx}
\Hat{u}^h_i := [I^h_H \Bar{u}^H_{\cdot}]_i,
\end{equation}
and the correction term can be approximated by omitting long-range($|j-i|>m$) terms,
\begin{equation} 
\label{eq:odd_corr_approx}
\Hat{e}^h_i := h \sum_{|j-i| \leq m} (G^{hh}_{i,j} - \Hat{G}^{hh}_{i,j}) f^h_j,
\end{equation}
then we have 
\begin{equation} 
\label{eq:odd_cor_multi-sum}
u^h_i \simeq \Hat{u}^h_i + \Hat{e}^h_i. \\
\end{equation}
Applying this idea recursively on a coarser grid till $\sqrt{n}$ points achieved, the total complexity will be reduced to $O(n\log n)$\cite{brandt1990multilevel}\cite{mlmi2005acoustics}. 

In summary, instead of calculating a dense matrix-vector product directly, MLMI approximates the product by first calculating it on a coarse grid and then using interpolations and corrections to refine it on finer grids. 


\subsection{GreenMGNet} \label{sec:gmg}
\begin{figure}[h]
\centering
\includegraphics[scale=0.12]{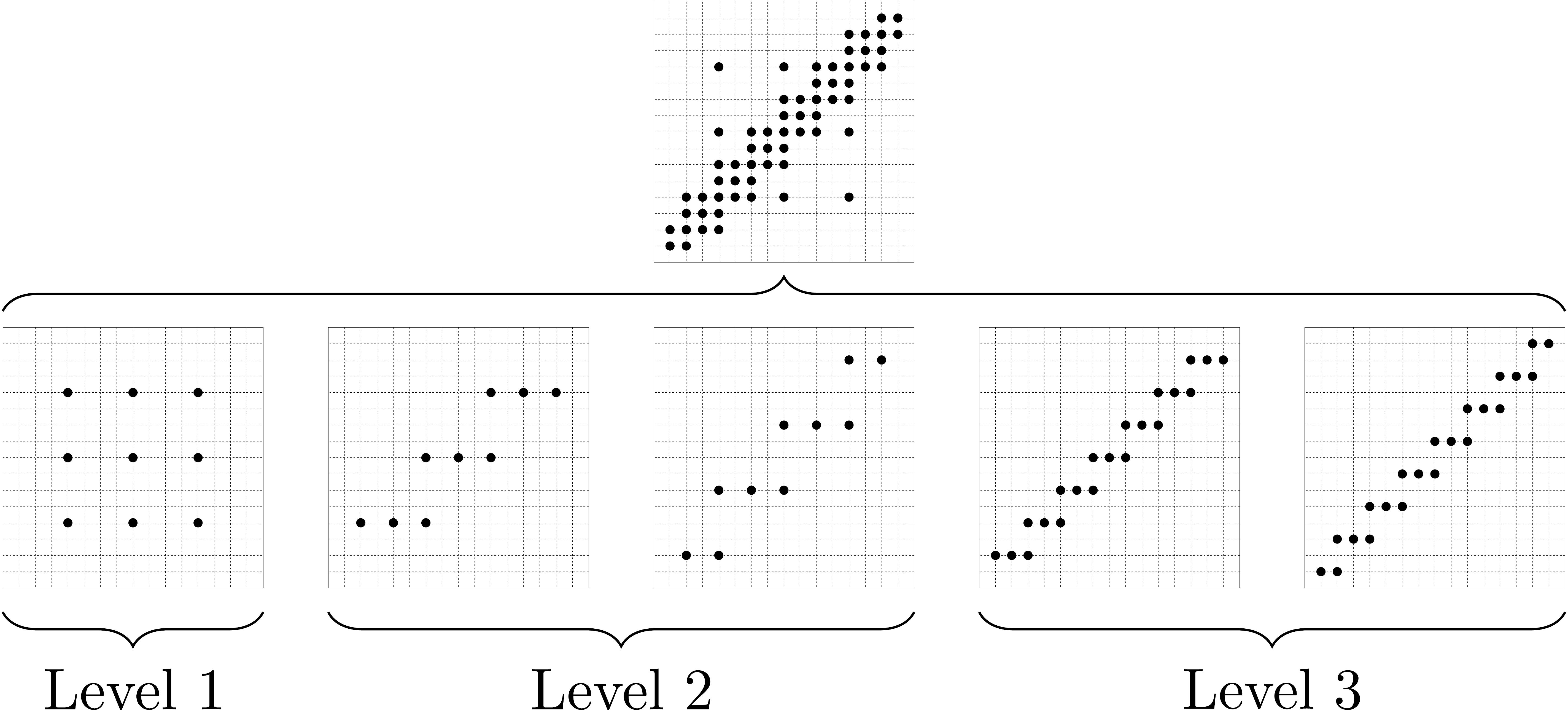}
\caption{MLMI only uses the subset of kernel values, which contains values on the coarsest grid points and the neighbors of diagonal points on each finer grid, for matrix-vector product calculation. Neighbors on each finer grid can be further classified as neighbors of $i=2I$, and neighbors of $i \neq 2I$, where $i$ and $I$ are relative indices between two levels. After decomposing the kernel values into subsets, it is possible to recover an approximant of the original full dense kernel by interpolation and corrections.}
\label{fig:mlmi_decomp}
\end{figure}

In the MLMI algorithm, only points near the diagonal or on a coarse grid are involved in matrix-vector product calculation(Figure \ref{fig:mlmi_decomp} gives an example of 3-level decomposition with local correction range $1$), which inspired us to leverage this idea in GL training to increase the training and inference efficiency under asymptotically smooth assumption. 

\begin{figure}[h!]
\centering
\includegraphics[scale=0.7]{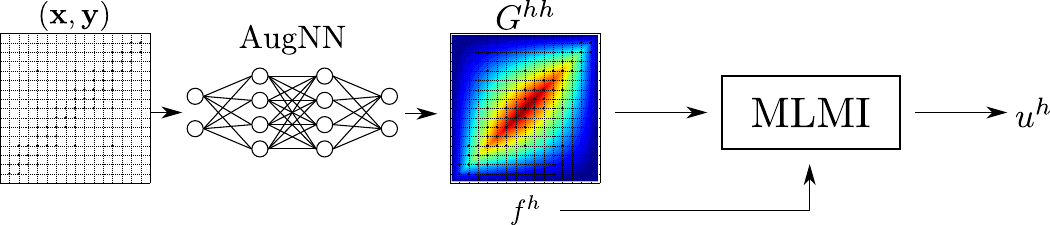}
\caption{Schematic of GreenMGNet}
\label{fig:gmgn}
\end{figure}

Figure \ref{fig:gmgn} shows the schematic of GreenMGNet. Similar to GL architecture, GreenMGNet can also be separated into a kernel evaluation step and a discrete kernel integration step. In the kernel evaluation step, only a subset of full grid points is needed. After determining the levels $k$ and local range $m$, the set of needed points is fixed and the corresponding kernel values are evaluated by AugNN. In the discrete kernel integration step, the solution is calculated by MLMI through multi-level interpolation and correction. During the training stage, if we coarsen the gridpoints proportional to $\sqrt{n}$, only $O(n\mathrm{log}n)$ points are needed, as a result, the operations of forward propagation of neural network and matrix-vector product will be reduced from $O(n^2)$ to $O(n\mathrm{log}n)$. Since the trained neural network is an approximated kernel function, both coarse kernel values and correction terms in (\ref{eq:even-cor_local})(\ref{eq:odd_corr_approx}) can be pre-computed. For an unseen force function $f$, we only need to apply MLMI to calculate solution $u$ with $O(n\mathrm{log}n)$ complexity.

\begin{figure}[h!]
  \centering
  \includegraphics[scale=0.8]{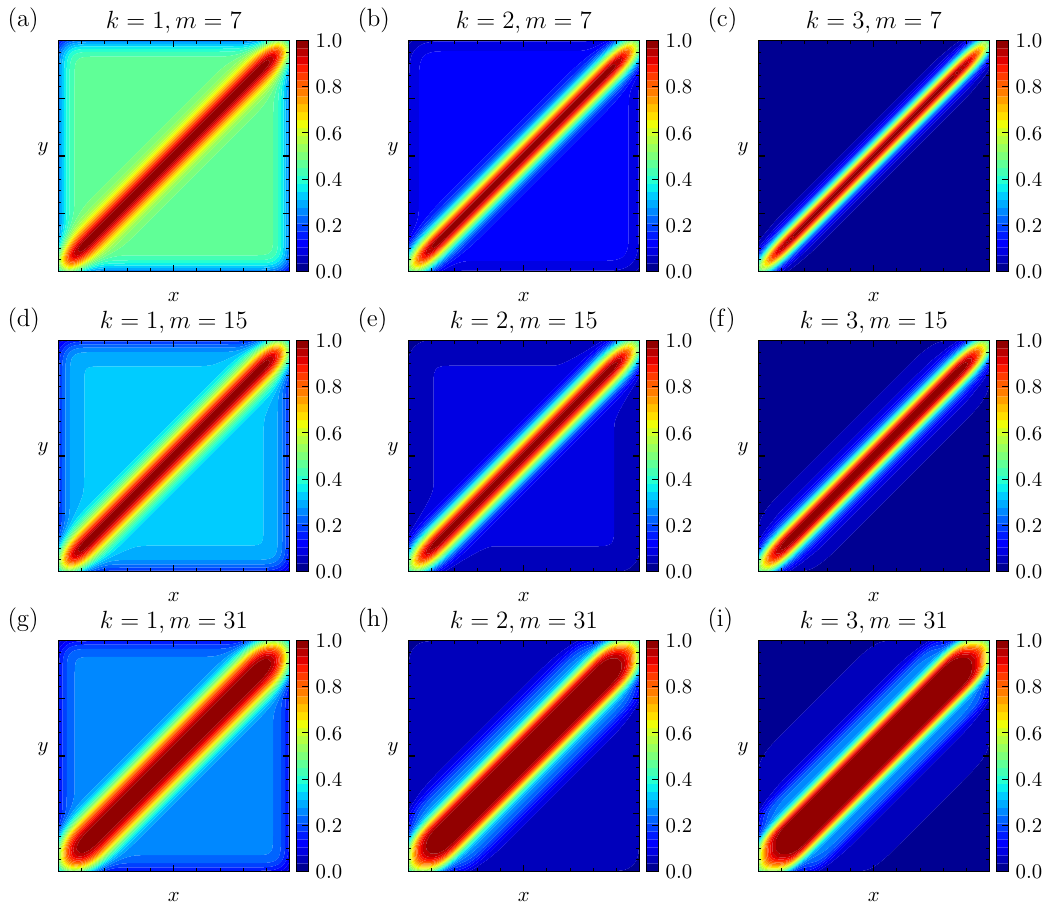}
  \caption{Normalized density functions of points used to train GreenMGNet. We approximate the probability density function by Gaussian kernels. The density function is normalized to $[0,1]$ for the convenience of visualization.}
  \label{fig:gmgn_distribution}
\end{figure}

From another viewpoint, GL learns a kernel function through point sampling from a uniform distribution. Nevertheless, GreenMGNet learns a kernel function by a set of points from a symmetric diagonally dominant distribution, determined by $m$ and $k$. Figure \ref{fig:gmgn_distribution} exhibits some normalized density functions of points used by MLMI. The density function reflects the significance of each location treated by GreenMGNet. More points near the diagonal are used in GreenMGNet, which is compatible with an asymptotically smooth assumption. It's important to mention that, despite some combinations of $m$ and $k$ having similar density functions, the amount of data points is different, which may lead to distinct results. 

\section{Numerical experiments} \label{sec:num_exps}

In this section, we present several numerical experiments to demonstrate the capability of GreenMGNet on various operator learning tasks under the asymptotic smoothness assumption. First, we exhibit the results of GL with Aug-NN(GL-aug). Then, we present the results of GreenMGNet, where Aug-NN is used as the default. For comparison, the GL serves as the baseline. In all experiments, the neural networks are fully connected neural networks(FNNs) with 4 hidden layers with 50 neurons in each layer. The activation function is set to Rational\cite{boulle2020rational} as GL. We train all the networks with relative $L_2$ error(\ref{eq:rl2}) by the Adam optimizer\cite{KingBa15}. For each case, five training trails with different random seeds are performed to compute the mean error and the standard deviation. All the computations are carried on a single Nvidia RTX A6000 GPU with 48GB memory. 

To evaluate the accuracy of different methods, we compute the relative $L_2$ error of solutions on the test sets and the relative $L_2$ error of kernels, if the exact kernel functions are available. Specifically, the relative $L_2$ error of solutions $\epsilon_u$ and the relative $L_2$ error of kernel functions $\epsilon_{G}$ are defined as,
\begin{equation} 
\label{eq:rl2_metric}
\epsilon_{u} = \frac{1}{N} \sum_{k=1}^{N} \frac{\| u^{(k)} - u_{\text{approx}}^{(k)} \|}{ \|u^{(k)} \| }, \epsilon_{G} = \frac{\| G - G_{\text{approx}} \|}{ \| G \| }, \\
\end{equation} where $u_{\text{approx}}$ and $G_{\text{approx}}$ denotes corresponding approximants. To visualize the results from different methods, we compute the pointwise absolute error of a sample solution or a kernel, denotes $E_u(\mathbf{x})$ and $E_G(\mathbf{x}, \mathbf{y})$ respectively, 
\begin{equation} 
\label{eq:abs}
E_{u}(\mathbf{x}) = | u^{(k)}(\mathbf{x}) - u_{\text{approx}}^{(k)}(\mathbf{x}) |, E_{G}(\mathbf{x}, \mathbf{y}) = | G(\mathbf{x}, \mathbf{y}) - G_{\text{approx}}(\mathbf{x}, \mathbf{y}) |. \\
\end{equation}

For one-dimensional tasks, the size of training and testing sets are both set to 100. We train the networks for 10,000 epochs with full training sets in each optimization step. A multi-step learning rate schedule is used. The initial learning rate is set to 0.01 and decays at 1,000th and 3,000th epochs with a factor of 0.1. The forcing terms of the one-dimensional tasks are drawn from a Gaussian process with a mean zero and a length scale of 0.03. The discrete forcing vectors are evaluated on a regular grid defined on $[0,1]$ with $n=513$. If the target problem is not defined on $[0,1]$, a simple linear transformation is applied.

For two-dimensional tasks, we train the networks for 1,000 epochs with an initial learning rate of 0.01, which decays at 100th and 300th epochs with a factor of 0.1. The sizes of the train and test sets are both set to 200. The 2D random forcing terms are generated by \texttt{randfun2} function of Chebfun software system\cite{driscoll2014chebfun, filip2019smooth} with a space scale of $0.2$. The 2D discrete forcing matrices are evaluated with a $65 \times 65$ cartesian grid defined on the unit square. Due to the limitations in computational resources, the batch size is set to $20$ in each optimization step and we randomly select up to $20\%$ points from the full grid points in either GL or GL-aug experiment during the training stage and use full grid points for validation. Unless otherwise specified, we obtain the solution by FEniCS software\cite{logg2012automated}.

To evaluate the efficiency of different methods, we measure the average execution time\footnote{\url{https://discuss.pytorch.org/t/how-to-measure-execution-time-in-pytorch/111458}} and peak GPU memory cost\footnote{\url{https://pytorch.org/docs/stable/generated/torch.cuda.max_memory_allocated.html}} for both training and inference. We generate a batch of 8 fake data with the same resolution as the real training dataset and run 1,000 iterations for 1D tasks and 100 iterations for 2D tasks. It is worth noting that the measured execution time during the inference stage doesn't include kernel evaluation, as the kernel should be pre-computed at the inference stage.

\subsection{Test problems}
GreenMGNet is designed to learn any unknown asymptotically smooth kernels, a typical application is to estimate solution operators of elliptic PDEs. The test problems are mainly from elliptic PDEs, containing several one-dimensional and two-dimensional tasks, and some have exact kernels. The details of the tasks are listed as follows,  
\begin{enumerate}
\item \textbf{One-dimensional Logarithm kernel}
The logarithm kernel is a typical potential-type kernel which is asymptotically smooth and has a singularity at the diagonal. The proposed method can be used to learn this kernel from data pairs calculated by the integral transform,
\begin{equation} 
\label{eq:log}
\begin{split}
u(x) & = \int_{-1}^{1} \mathrm{ln}|y-x| f(x) \mathrm{d}y, \quad x \in [-1,1]. \\
\end{split}
\end{equation}
The corresponding discrete kernel can be calculated analytically, given as, 
\begin{equation} 
\label{eq:log-coef}
\begin{split}
G^{hh}_{i,j} &= \frac{1}{h} \int_{y_j-\frac{h}{2}}^{y_j+\frac{h}{2}} \mathrm{ln}|x_i-y| \mathrm{d}y \\
&= \frac{1}{h} \Bigl( (x_i-y_j+\frac{h}{2})(\mathrm{ln}|x_i-y_j+\frac{h}{2}|-1) \\
&\quad - (x_i-y_j-\frac{h}{2})(\mathrm{ln}|x_i-y_j-\frac{h}{2}|-1) \Bigl).
\end{split}
\end{equation}
We obtained a discrete approximation of $u(x)$ by directly calculating the product of $G^{hh}$ and discrete forcing functions.
\item \textbf{One-dimensional Poisson equation} Poisson equation is a uniformly elliptic PDE. The corresponding Green's function is known and has a cusp on the diagonal. We have discussed this equation in Section \ref{para:poisson_1d}, which is given as
\begin{equation} 
\label{eq:poisson1d_}
\begin{split}
-\frac{d^2 u}{d x^2} & = f(x), \quad x \in [0,1] \\
u(0) & = u(1) = 0.
\end{split}
\end{equation}    
\item \textbf{One-dimensional Schr\"{o}dinger equation} We consider a steady-state one-dimensional Schr\"{o}dinger equation with double-well potential $V(x) = x^2 + 1.5 \exp{(-(4x)^4)}$, which is an elliptic PDE, 
\begin{equation} 
\label{eq:schodinger}
\begin{split}
-h^2 \frac{d^2 u}{ d x^2} + V(x)u &= f, \quad\quad x \in [-3, 3] \\
u(-3) = u(3) &= 0,
\end{split}
\end{equation}
where $h=0.1$ is a coefficient of the equation not meshsize.

\item \textbf{One-dimensional Airy's equation} The last one-dimensional task is Airy's equation with homogeneous Dirichlet boundary conditions and with parameter $\theta = 10$. This task is also an elliptic PDE operator learning task:
    \begin{equation} 
    \label{eq:airy}
    \begin{split}
    -\frac{d^2 u}{d x^2} + \theta^2 x u &= f, \quad\quad x \in [0, 1] \\
    u(0) = u(1) &= 0,
    \end{split}
    \end{equation}
    \item \textbf{Two-dimensional Poisson equation}
    We consider a two-dimensional Poisson problem defined on a unit square,     
    \begin{equation} 
    \label{eq:poissonrect}
    \begin{split}
    -\nabla \cdot (\nabla u(\mathbf{x})) &= f(\mathbf{x}), \quad \mathbf{x} \in \Omega, \\
    u(\mathbf{x}) &= 0, \quad \mathbf{x} \in \partial \Omega,
    \end{split}
    \end{equation}
    where $\Omega=[0,1] \times [0,1]$. Due to $\Omega$ being a unit square, no analytical Green's function is available in this case, which is different from what we have discussed in Section \ref{para:poisson_2d}.  
    \item \textbf{Two dimensional Darcy's Flow}
    The last example is Darcy's Flow which is a typical uniformly elliptic PDE,
    \begin{equation} 
    \label{eq:darcyrect}
    \begin{split}
    -\nabla \cdot (A(\mathbf{x}) \nabla u(\mathbf{x})) & = f(\mathbf{x}), \quad \mathbf{x} \in \Omega \\
    u(\mathbf{x}) & = 0, \quad\quad \mathbf{x}\in \partial \Omega.
    \end{split}
    \end{equation}
    where $\Omega=[0,1] \times [0,1]$. The coefficient $A(\mathbf{x})$ is defined as $a(\mathbf{x}) I$, where $I$ is a $2 \times 2$ identity matrix and $a(\mathbf{x})$ is a random scalar output function sampled from a Gaussian random field followed by an exponential transform.
\end{enumerate}

\subsection{Accuracy of AugNN} \label{subsec:aug_nn-result}

\begin{figure}[h!]
    \centering
    \includegraphics[scale=0.7]{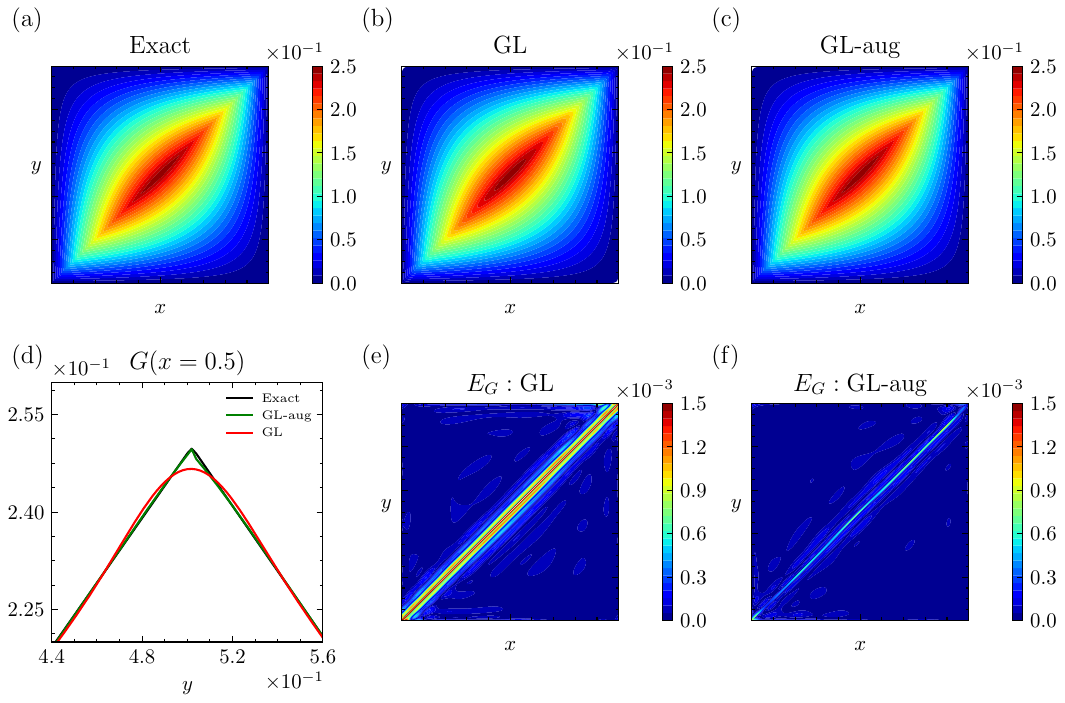}
    \caption{Result of Poisson 1D problem: (a) Exact Green's function. (b) Green's function learned by GL method. (c) Green's function learned by GL-aug. (d) The slice $G(x=0.5)$ comparison. (e) Absolute error $E_G$ of GL. (e) Absolute error $E_G$ of GL-aug.}
    \label{fig:poisson1d_result}
\end{figure}

\begin{figure}[h!]
    \centering
    \includegraphics[scale=0.7]{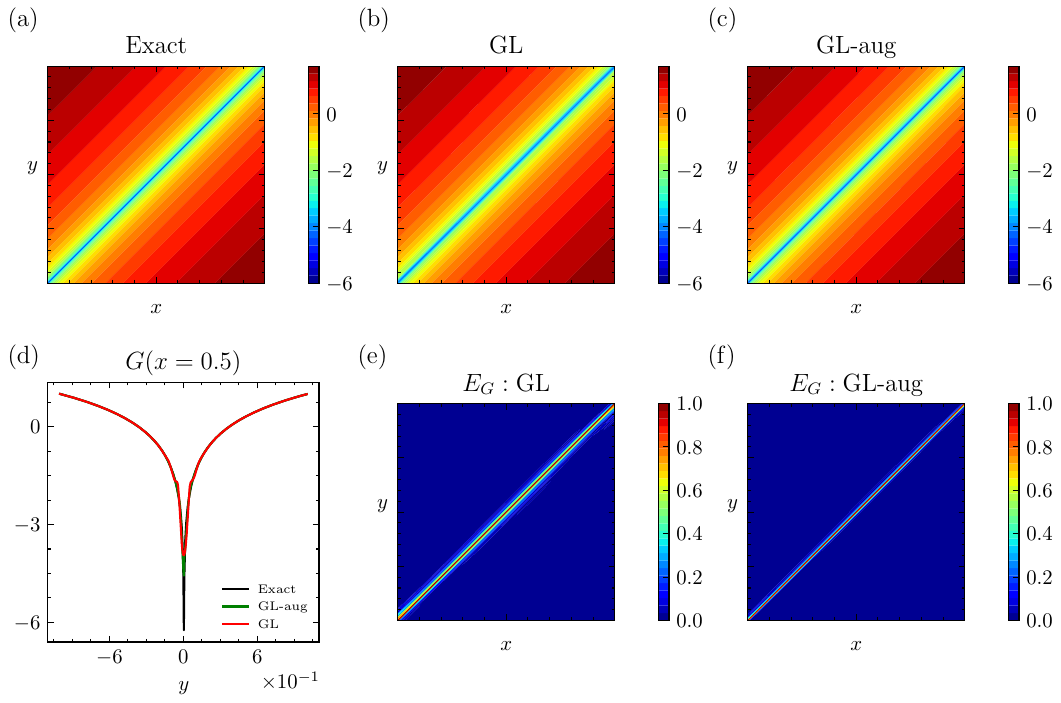}
    \caption{Result of Logarithm 1D problem: (a) Exact Green's function. (b) Green's function learned by GL method. (c) Green's function learned by GL-aug. (d) The slice $G(x=0.5)$ comparison. (e) Absolute error $E_G$ of GL. (e) Absolute error $E_G$ of GL-aug.}
    \label{fig:log1d_result}
\end{figure}

We first consider the problems with exact kernels. Table \ref{table:analytic1d_aug} shows the results of both $\epsilon_u$ and $\epsilon_G$. Compared to the baseline, GL-aug achieves a remarkable improvement of $26.56 \%$ and $71.21 \%$ on $\epsilon_u$ and $\epsilon_G$ for the Poisson 1D problem, and a $20.03 \%$ and $26.67 \%$ improvement on $\epsilon_u$ and $\epsilon_G$ for the Logarithm 1D problem. Learning the logarithm kernel is difficult due to the singularities on the diagonal. While GL-aug exhibits significant improvement, discretization and optimization error remain a major source of error as mentioned in \cite{boulle2023ellipticefficient}. In contrast, Green's function of the Poisson 1D problem has a cusp on the diagonal, which is easier to learn. In both cases, the improvement stems from a better approximation near the diagonal. A visual analysis is depicted in Figure \ref{fig:poisson1d_result} and Figure \ref{fig:log1d_result}. As one can see, the GL-aug method has a smaller $E_G$ near the diagonal. GL method tends to learn an everywhere smooth approximation, but GL-aug can preserve discontinuity which leads to a better result.

\begin{table}[h!]
\centering
\begin{tabular}{ lcccc }
    \toprule
         & \multicolumn{2}{c}{GL}    & \multicolumn{2}{c}{GL-aug}\\
    \cmidrule(lr){2-3}\cmidrule(lr){4-5}
         & $\epsilon_G$ & $\epsilon_u$ & $\epsilon_G$ & $\epsilon_u$ \\
    \midrule
    Poisson & $ 0.271 \pm 0.008 $ & $0.064 \pm 0.006$ & $ \mathbf{0.078 \pm 0.012} $ & $ \mathbf{0.047 \pm 0.003} $ \\
    Logarithm & $ 10.62 \pm 0.285 $ & $ 0.075 \pm 0.005 $ & $ \mathbf{8.492 \pm 0.538} $ & $ \mathbf{0.055 \pm 0.006} $ \\
    \bottomrule
\end{tabular}
\caption{Relative errors comparison($\times 10^{-2}$) for problems with exact kernels}
\label{table:analytic1d_aug}
\end{table}

\begin{figure}[h!]
    \centering
    \includegraphics[scale=0.5]{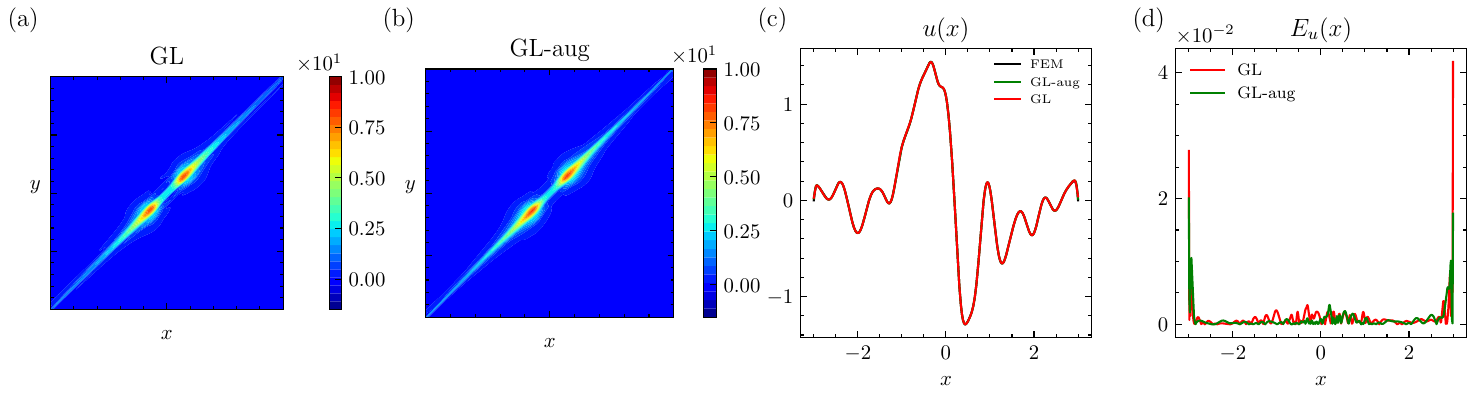}
    \caption{Result of Schr\"{o}dinger 1D problem: (a) Green's function learned by GL method. (b) Green's function learned by GL-aug method. (c) Sample solution comparison. (d) Absolute error $E_u$ comparison.}
    \label{fig:schrodinger1d_result}
\end{figure}

\begin{figure}[h!]
    \centering
    \includegraphics[scale=0.5]{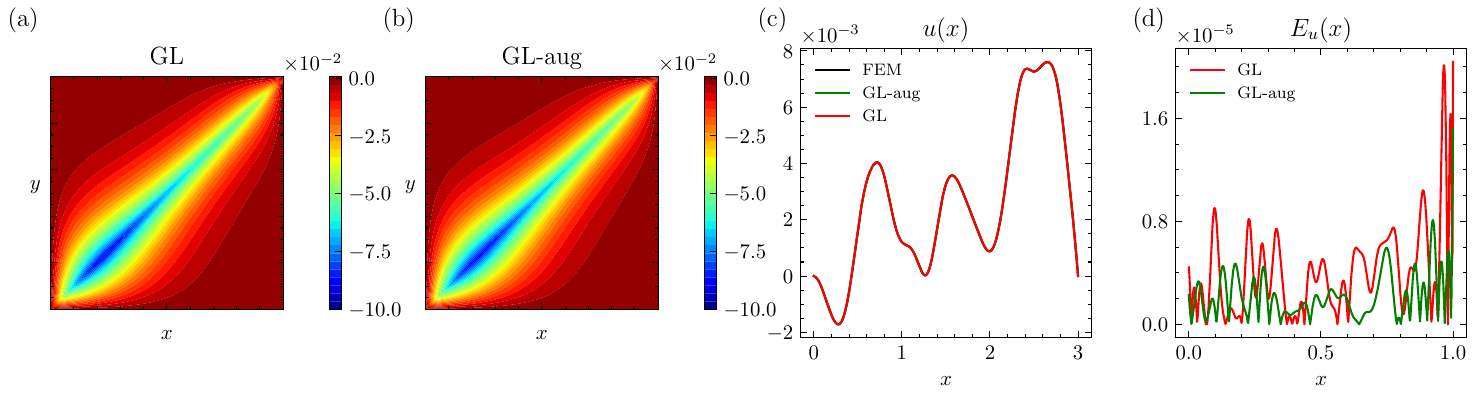}
    \caption{Result of Airy 1D problem: (a) Green's function learned by GL method. (b) Green's function learned by GL-aug method. (c) Sample solution comparison. (d) Absolute error $E_u$ comparison.}
    \label{fig:airy1d_result}
\end{figure}

For problems without exact kernels, the results are summarized in Table \ref{table:other_aug}. GL-aug achieves better results in terms of $\epsilon_u$ in all cases. Visualizing the learned kernels, we observe different types of non-smoothness on the diagonals. The kernel function of 1D Schr\"{o}dinger equation has two poles on the diagonal, see Figure \ref{fig:schrodinger1d_result}. Both GL and GL-aug can capture the poles; however, GL-aug provides a sharper approximation around the poles, outperforming GL $25.92 \%$ on $\epsilon_u$. The learned kernel of 1D Airy's equation has cusps on the diagonal, similar to the Green's function of the 1D Poisson equation, but the location and sharpness of cusps on the diagonal are different, see Figure \ref{fig:airy1d_result}. For 1D Airy's equation, the result of GL-aug is $7.79 \%$ better than GL on $\epsilon_u$. For the 2D Poisson equation and Darcy's Flow on the unit square, GL-aug can learn kernel function with sharper diagonals. Compared to 1D problems, the kernel function of a 2D problem has a lower proportion of points on the diagonal. Therefore, the improvement is not very significant but amounts to $4.50 \%$  and $1.41 \%$ improvement on the Poisson equation and Darcy's flow respectively.

In summary, GL tends to learn a smooth kernel function, even non-smoothness exists near the diagonal. In contrast, GL-aug can capture the non-smoothness more accurately on the diagonal by giving a sharp approximation. Although only part of the diagonal contains non-smoothness in some cases, GL-aug still facilitates the result. 

\begin{figure}[h!]
    \centering
    \includegraphics[scale=0.7]{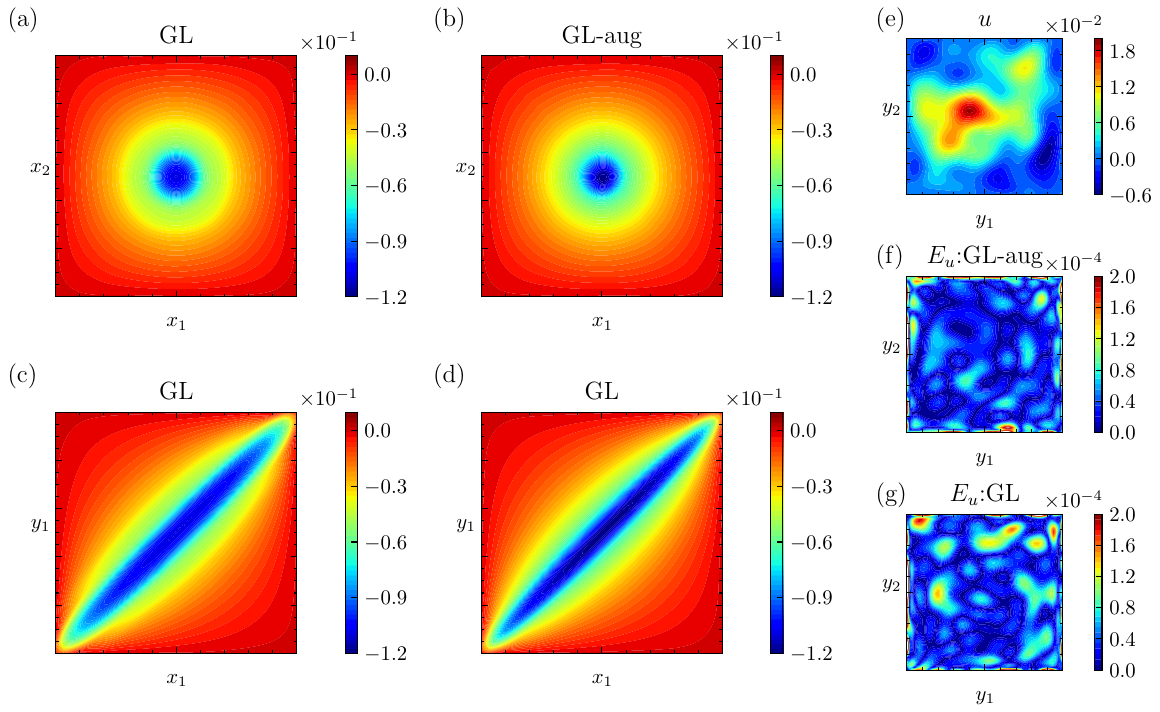}
    \caption{Result of Poisson 2D problem: (a) $G(x_1, x_2, 0, 0)$ learned by GL method. (b) $G(x_1, x_2, 0, 0)$ learned by GL-aug method. (c) $G(x_1, 0, y_1, 0)$ learned by GL method. (d) $G(x_1, 0, y_1, 0)$ learned by GL-aug method. (e) Sample solution $u$. (f) Absolute error $E_u$ of GL-aug method. (g) Absolute error $E_u$ of GL method.}
    \label{fig:poissonrect_result}
\end{figure}

\begin{figure}[h!]
    \centering
    \includegraphics[scale=0.7]{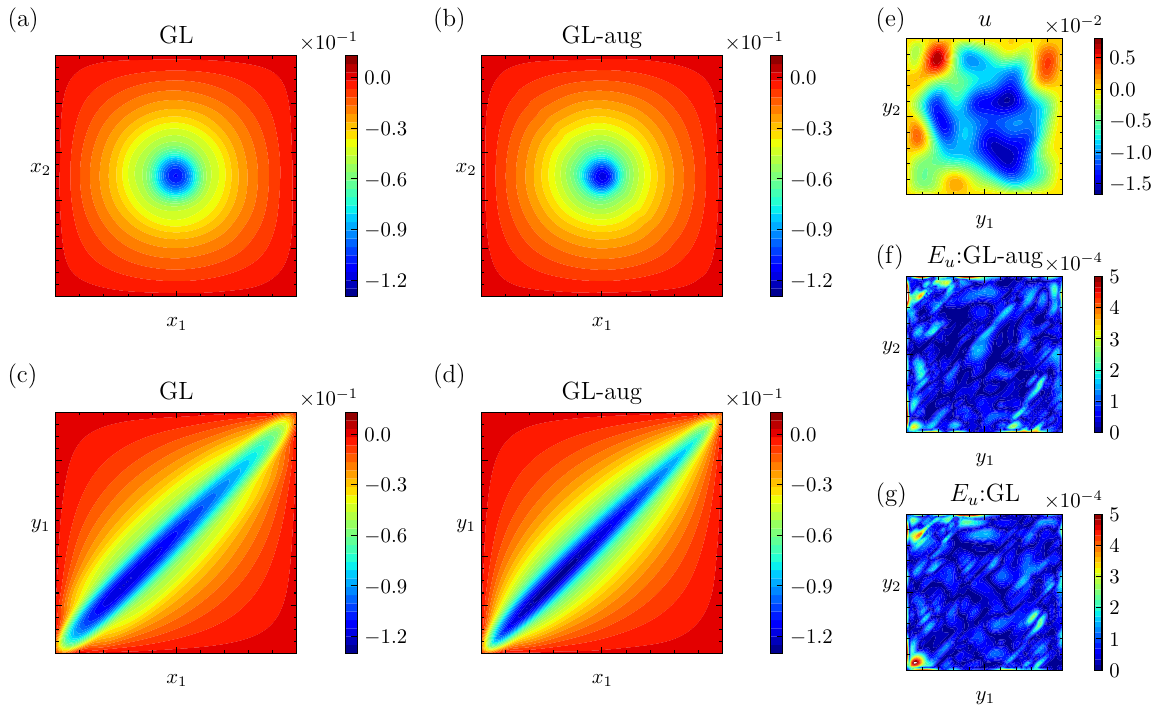}
    \caption{Result of Darcy's flow 2D problem: (a) $G(x_1, x_2, 0, 0)$ learned by GL method. (b) $G(x_1, x_2, 0, 0)$ learned by GL-aug method. (c) $G(x_1, 0, y_1, 0)$ learned by GL method. (d) $G(x_1, 0, y_1, 0)$ learned by GL-aug method. (e) Sample solution $u$. (f) Absolute error $E_u$ of GL-aug method. (g) Absolute error $E_u$ of GL method.}
    \label{fig:darcyrect_result}
\end{figure}

\begin{table}[h!]
\centering
\begin{tabular}{ccc}
\toprule 
    & GL    & GL-aug \\
\midrule
    Schr\"{o}dinger & $ 0.189 \pm 0.020 $ & $\mathbf{0.140 \pm 0.034}$ \\
    Airy & $ 0.077 \pm 0.012 $ & $ \mathbf{0.071 \pm 0.008} $ \\
    Poisson 2D & $ 0.867 \pm 0.170 $ & $ \mathbf{0.828 \pm 0.096} $ \\
    Darcy's Flow & $ 1.202 \pm 0.107 $ & $ \mathbf{1.185 \pm 0.089} $ \\
\bottomrule
\end{tabular}
\caption{$\epsilon_u$ ($\times 10^{-2}$) comparison for problems without analytical kernels}
\label{table:other_aug}
\end{table}

\subsection{Accuracy of GreenMGNet}

Under the asymptotically smooth assumption of the underlying kernel, GreenMGNet leverages AugNN and MLMI in GL training. The local correction range $m$ and the number of levels $k$ are two key parameters in GreenMGNet, which determine the distribution and the amount of points used during training as we discussed in Section \ref{sec:gmg}. For 1D problems, we choose $m$ from the set $\{0, 1, 3, 7, 15, 31\}$, and for 2D problems, $m$ is selected from the set $\{0, 1, 3, 5\}$, where $m=0$ indicates that no local correction involved. The coarse level $k$ is limited to a maximum of 3. Table \ref{table:gmgn_main} shows the comparison of GL, GL-aug, and GreenMGNet(best result from $(m,k)$ pairs). GreenMGNet surpasses GL in all cases and achieves an even better $\epsilon_u$ in most cases compared to GL-aug. The minimum accuracy improvement relative to GL is 2D Darcy's Flow, which shows a gain of $3.8\%$. The maximum accuracy improvement relative to GL occurs in 1D Smeanchr\"odinger, with an impressive $39.15\%$ improvement.

\begin{table}[h!]
\centering
\begin{tabular}{cccccc}
\toprule 
    & $k$ & $m$ & GreenMGNet & GL & GL-aug \\
\midrule
    Poisson 1D  & 1 & 3 & $\mathbf{0.046 \pm 0.003}$ & $0.064 \pm 0.006$ & $ 0.047 \pm 0.003 $  \\
    Schr\"{o}dinger & 3 & 31 & $\mathbf{0.115 \pm 0.017}$ & $ 0.189 \pm 0.020 $ & $0.140 \pm 0.034$\\
    Airy & 2 & 15 & $0.074\pm0.008$ & $ 0.077 \pm 0.012 $ & $\mathbf{0.071 \pm 0.008}$ \\
    Poisson 2D & 1 & 5 & $\mathbf{0.780\pm0.120}$ & $ 0.867 \pm 0.170 $ & $0.828 \pm 0.096$ \\
    Darcy's Flow & 1 & 5 & $\mathbf{1.156\pm0.120}$ & $ 1.202 \pm 0.107 $ & $1.185 \pm 0.089$\\
\bottomrule
\end{tabular}
\caption{$\epsilon_u$ ($\times 10^{-2}$) comparison for the test problems; the first two columns are the corresponding parameters of the best GreenMGNet setting.}
\label{table:gmgn_main}
\end{table}

We observe that the best GreenMGNet configuration for different tasks requires different parameter settings. This indicates that the amount and distribution of data points used to learn kernels should vary according to the task. Since GreenMGNet uses fewer points, we compare it with both GL and GL-aug, which are trained with $p$ percentage of points randomly sampled from full grid points. 
The percentage $p$ is chosen from set $\{0.03, 0.05, 0.07, 0.1, 0.15, 0.25, 0.3, 0.5, 1.0\}$ for 1D tasks and $\{ 0.001, 0.005, 0.03, 0.05, 0.1, 0.15, 0.2 \}$ for 2D tasks. For GreenMGNet, $p$ is determined by $(m,k)$ pair; therefore, the data points are selected not randomly but according to the MLMI.

\begin{figure}[h!]
    \centering
    \includegraphics[scale=0.7]{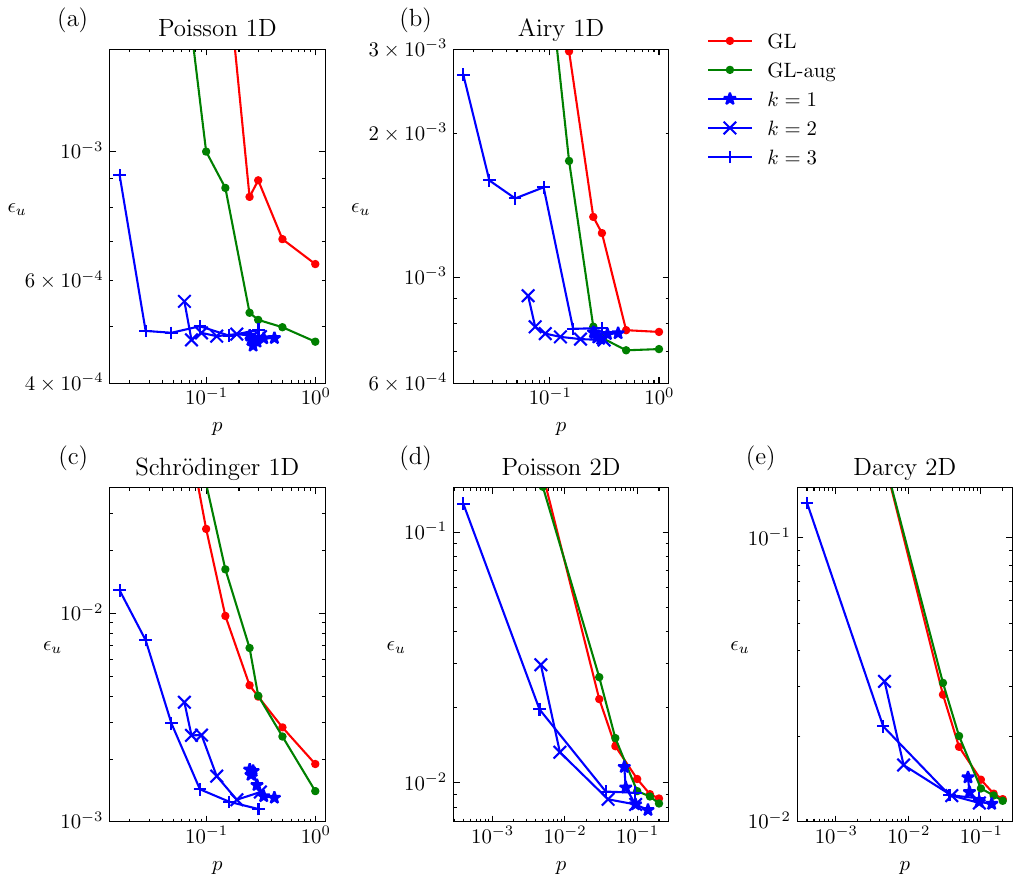}
    \caption{Relationship of $\epsilon_u$ and $p$, where $p$ is the percentage of points. For GreenMGNet, we can calculate $p$ after setting $(k, m)$ by counting the amounts of points used by MLMI.}
    \label{fig:gmgn_result}
\end{figure}

Figure \ref{fig:gmgn_result} displays the relationship between $\epsilon_u$ and $p$ for different methods. The $\epsilon_u$ of GreenMGNet converges faster than that of both GL and GL-aug as $p$ increases. This trend verifies the underlying asymptotically smooth assumption, showing that training points from multi-level diagonal dominant structures will improve learning performance compared to uniform random sampling.

For 1D Poisson and 1D Airy problem, there are cusps on the diagonal of the corresponding kernel function, they are not as sharp as singularities, see Figure \ref{fig:poisson1d_result} and Figure \ref{fig:airy1d_result}. In this case, a small range of correction $m$ is enough to learn a good approximation, and increasing the range of correction does not result in significant improvement. 

In contrast, for 1D Schr\"{o}dinger, 2D Poisson, and 2D Darcy's Flow, singularities are more likely exists on the diagonal, see Figure \ref{fig:schrodinger1d_result}, Figure \ref{fig:poissonrect_result} and Figure \ref{fig:darcyrect_result}, and thus the trend curve does not converge as quickly as 1D Poisson problem. However, it significantly outperforms both GL and GL-aug. We observe that, despite different coarse levels, increasing the correction range $m$ will achieve better results. We also notice the amount of data points still matter; fewer coarse levels mean more data points involved, and consequently, the results are better.

\subsection{Efficiency of GreenMGNet}

We present the representative results from \ref{appendix_full_results} in Table \ref{table:sub_efficiency_1d} and Table \ref{table:sub_efficiency_2d}. GL-aug uses a slightly larger neural network than GL due to the augmented output, which typically results in slightly higher memory costs and longer training time. However, in practice, due to the implementation of neural networks by deep learning frameworks, the measured allocated memory of GL-aug and GL could be the same because the number of parameters of the two models is almost identical.

GreenMGNet demonstrates remarkable efficiency compared to both GL and GL-aug. Specifically, for 1D tasks, GreenMGNet($k=2, m=7$) achieves overall better accuracy than GL using only $12\%$ of the full grid data. The training time is reduced to $54.10\%$, and training GPU memory cost is only $7.5\%$. The inference GPU memory cost is further reduced to $10.0\%$. For the 2D Poisson problem, GreenMGNet($k=2,m=5$) achieves better accuracy than both GL and GL-aug using only $9.5\%$ of full grid data. The training time is shortened to $62.3\%$, and training GPU memory cost is only $37.5\%$. The inference GPU memory cost is reduced to $42.8\%$.

During the training stage, execution time and memory costs primarily arise from kernel evaluation. Consequently, GreenMGNet exhibits significant improvements compared to GL and GL-aug because of fewer points requirement. During the inference stage, discrete kernel integration is the main consumption. Unlike direct matrix-vector multiplication, the implementation of MLMI is not highly optimized for GPU devices, therefore, inference execution time is slower than GL and GL-aug. However, the memory cost of GreenMGNet remains significantly lower than GL and GL-aug. Notably, in 2D cases, evaluating the full kernel directly is not feasible due to GPU memory limitations. Therefore, we split the kernel into several kernel blocks and then evaluate the integral of each part sequentially for GL and GL-aug. In contrast, such splitting is unnecessary for GreenMGNet. 

\begin{table}[h!]
\begin{adjustbox}{width={\textwidth},totalheight={\textheight},keepaspectratio}
\centering
\begin{tabular}{ccccccccccc}
\toprule 
         & & & & \multicolumn{2}{c}{Training}    & \multicolumn{2}{c}{Inference} & \multicolumn{3}{c}{$\epsilon_u$($\times 10^{-2}$)}\\
    \cmidrule(lr){5-6}\cmidrule(lr){7-8} \cmidrule(lr){9-11}
           & $k$ & $m$ & $p$ & Time(ms)  & Memory(Mb) & Time(ms)  & Memory(Mb) & Poisson 1D & Airy 1D & Schr\"odinger 1D \\
\midrule
GreenMGNet & 2 & 7 & 0.12 & 28.78 & 248.87 & 2.19 & 217.26 & $0.048\pm0.004$ & $0.075\pm0.010$ & $0.166\pm0.036$ \\
\midrule
GL      & - & - & 1.00 & 53.20 & 2483.21 & 0.05 & 1621.72 & $0.064\pm0.006$ & $0.077\pm0.012$ & $0.189\pm0.020$ \\
\midrule
GL-aug  & - & - & 1.00 & 54.39 & 2483.21 & 0.05 & 1621.72 &  $0.047\pm0.003$ & $0.071\pm0.008$ & $0.140\pm0.034$\\
\bottomrule
\end{tabular}
\end{adjustbox}
\caption{GreenMGNet results comparison of the models for 1D tasks.}
\label{table:sub_efficiency_1d}
\end{table}

\begin{table}[h!]
\begin{adjustbox}{width={\textwidth},totalheight={\textheight},keepaspectratio}
\centering
\begin{tabular}{cccccccccc}
\toprule 
         & & & & \multicolumn{2}{c}{Training}    & \multicolumn{2}{c}{Inference} & \multicolumn{2}{c}{$\epsilon_u$($\times 10^{-2}$)} \\
    \cmidrule(lr){5-6}\cmidrule(lr){7-8}\cmidrule(lr){9-10}         
           & $k$ & $m$ & $p$ & Time(ms)  & Memory(Mb) & Time(ms)  & Memory(Mb) & Poisson 2D & Darcy 2D \\
\midrule
GreenMGNet & 2 & 5 & 0.0955 & 417.21 & 12845.02 & 7.55 & 11676.11 & $0.820\pm0.178$ & $1.165\pm0.123$ \\
\midrule
GL  & - & - & 0.2000 & 668.88 & 34270.73 & 0.21 & 27304.10 &$0.867\pm0.170$ & $1.202\pm0.107$ \\
\midrule 
GL-aug & - & - & 0.2000 & 676.43 & 34270.73 & 0.21 & 27362.86 & 
 $0.828\pm0.096$ & $1.185\pm0.089$ \\
\bottomrule
\end{tabular}
\end{adjustbox}
\caption{Numerical results comparison of the models for 2D tasks.}
\label{table:sub_efficiency_2d}
\end{table}

\section{Conclusion} \label{sec:conclusion}

In this paper, we introduce GreenMGNet, an approach for approximating unknown asymptotically smooth kernels. The crucial idea is to model the asymptotically smooth kernel as a piecewise function using an augmented output neural network(AugNN) to preserve non-smoothness. To further enhance efficiency, we employ the Multi-Level Multi-Integration(MLMI) algorithm for discrete kernel evaluation. Our results demonstrate that GreenMGNet is both efficient and accurate for asymptotically smooth kernel estimation, making it suitable for learning solution operators of PDEs. 

Like many previous works, the errors of the present work are mainly from numerical integration and optimization, preventing us from achieving a theoretical order of accuracy. As a forthcoming extension, we plan to explore replacing the deep neural network with a shallow neural network and train it using greedy training algorithms, which can be rigorously proven to converge.

\section{Acknowledgments}
Y.J. Lee was funded in part by NSF-DMS 2208499.


\bibliographystyle{cas-model2-names}

\bibliography{GreenMGNet}

\appendix
\section{Appendix}

\subsection{Full Numerical Results} \label{appendix_full_results}

\begin{table}[h!]
\begin{adjustbox}{width={\textwidth},totalheight={\textheight},keepaspectratio}
\centering
\begin{tabular}{ccccccccccc}
\toprule 
         & & & & \multicolumn{2}{c}{Training}    & \multicolumn{2}{c}{Inference} & \multicolumn{3}{c}{$\epsilon_u$($\times 10^{-2}$)}\\
    \cmidrule(lr){5-6}\cmidrule(lr){7-8} \cmidrule(lr){9-11}
           & $k$ & $m$ & $p$ & Time(ms)  & Memory(Mb) & Time(ms)  & Memory(Mb) & Poisson 1D & Airy 1D & Schr\"odinger 1D \\
\midrule
GreenMGNet & 1 & 0  & 0.25 & 17.23 & 637.90 & 0.20 & 417.70 &  $0.048\pm0.004$ & $0.077\pm0.009$ & $0.178\pm0.020$ \\
           & 1 & 1  & 0.26 & 35.13 & 638.03 & 1.18 & 428.52 & $0.047\pm0.003$ & $0.076\pm0.011$ & $0.168\pm0.024$ \\
           & 1 & 3  & 0.27 & 24.63 & 638.21 & 1.19 & 448.80 & $0.046\pm0.003$ & $0.075\pm0.009$ & $0.175\pm0.027$ \\
           & 1 & 7  & 0.29 & 24.58 & 638.57 & 1.21 & 487.98 & $0.047\pm0.003$ & $0.074\pm0.010$ & $0.150\pm0.033$ \\
           & 1 & 15 & 0.34 & 31.07 & 639.27 & 1.31 & 560.36 & $0.048\pm0.003$ & $0.076\pm0.011$ & $0.133\pm0.037$ \\
           & 1 & 31 & 0.42 & 29.22 & 803.92 & 1.31 & 700.06 & $0.048\pm0.002$ & $0.076\pm0.009$ & $0.130\pm0.041$ \\
\midrule
           & 2 & 0 & 0.06 & 7.46 & 175.21 & 0.34 & 114.19 & $0.055\pm0.003$ & $0.091\pm0.009$ & $0.374\pm0.009$ \\
           & 2 & 1 & 0.07 & 23.12 & 175.38 & 2.12 & 130.71 & $0.047\pm0.003$ & $0.079\pm0.010$ & $0.260\pm0.017$ \\
           & 2 & 3 & 0.09 & 30.33 & 178.16 & 2.13 & 159.72 &  $0.049\pm0.003$ & $0.076\pm0.012$ & $0.261\pm0.016$  \\
           & 2 & 7 & 0.12 & 28.78 & 248.87 & 2.19 & 217.26 & $0.048\pm0.004$ & $0.075\pm0.010$ & $0.166\pm0.036$ \\
           & 2 & 15 & 0.19 & 26.77 & 383.69 & 2.31 & 325.94 & $0.049\pm0.002$ & $0.074\pm0.008$ & $0.128\pm0.024$  \\
           & 2 & 31 & 0.31 & 30.43 & 628.66 & 2.36 & 525.56 & $0.048\pm0.002$ & $0.074\pm0.009$ & $0.139\pm0.036$ \\
\midrule
           & 3 & 0 & 0.02 & 10.66 & 57.24 & 0.38 & 35.74 & $0.091\pm0.003$ & $0.265\pm0.213$ & $1.291\pm0.002$ \\
           & 3 & 1 & 0.03 & 43.14 & 66.37 & 3.05 & 55.06 & $0.049\pm0.003$ & $0.159\pm0.219$ & $0.744\pm0.005$ \\
           & 3 & 3 & 0.05 & 35.84 & 106.67 & 3.09 & 88.78 & $0.049\pm0.003$ & $0.146\pm0.190$ & $0.296\pm0.006$ \\
           & 3 & 7 & 0.09 & 31.58 & 185.63 & 3.16 & 154.49 &  $0.050\pm0.004$ & $0.154\pm0.227$ & $0.143\pm0.018$ \\
           & 3 & 15 & 0.16 & 34.88 & 336.77 & 3.20 & 277.30 &  $0.048\pm0.002$ & $0.078\pm0.011$ & $0.125\pm0.023$ \\
           & 3 & 31 & 0.30 & 37.73 & 607.49 & 3.36 & 504.39 &  $0.049\pm0.004$ & $0.078\pm0.009$ & $0.115\pm0.017$  \\
\midrule
GL & - & - & 0.03 & 6.95  & 96.17   & 0.05 & 1621.72 & 
   $10.49\pm6.578$ & $19.56\pm18.75$ & $50.15\pm55.83$ \\
   & - & - & 0.05 & 4.60  & 145.08  & 0.05 & 1621.72 & $2.559\pm2.678$ & $3.386\pm1.752$ & $10.28\pm2.876$ \\
   & - & - & 0.07 & 5.35  & 192.75  & 0.05 & 1621.72 & $0.869\pm0.930$ & $1.202\pm0.903$ & $6.889\pm5.521$ \\
   & - & - & 0.10 & 6.69  & 268.50  & 0.05 & 1621.72 & $0.240\pm0.221$ & $0.521\pm0.566$ & $2.542\pm2.498$ \\
   & - & - & 0.15 & 10.21 & 389.52  & 0.05 & 1621.72 & $0.221\pm0.178$ & $0.297\pm0.321$ & $0.972\pm0.430$ \\ 
   & - & - & 0.25 & 14.56 & 639.69  & 0.05 & 1621.72 & $0.084\pm0.022$ & $0.134\pm0.107$ & $0.452\pm0.123$ \\
   & - & - & 0.30 & 17.11 & 758.09  & 0.05 & 1621.72 & $0.089\pm0.039$ & $0.124\pm0.087$ & $0.399\pm0.135$ \\
   & - & - & 0.50 & 27.79 & 1257.95 & 0.05 & 1621.72 & $0.071\pm0.007$ & $0.077\pm0.012$ & $0.284\pm0.061$ \\
   & - & - & 1.00 & 53.20 & 2483.21 & 0.05 & 1621.72 & $0.064\pm0.006$ & $0.077\pm0.012$ & $0.189\pm0.020$ \\
\midrule
GL-aug & - & - & 0.03 & 4.30  & 96.17   & 0.05 & 1621.72 & $3.255\pm1.312$ & $4.496\pm1.873$ & $38.17\pm9.916$ \\
       & - & - & 0.05 & 4.95  & 145.08  & 0.05 & 1621.72 &  $0.511\pm0.556$ & $1.640\pm1.493$ & $15.59\pm7.209$ \\
       & - & - & 0.07 & 8.22  & 192.75  & 0.05 & 1621.72 &       $0.174\pm0.086$ & $0.568\pm0.510$ & $8.828\pm5.427$ \\
       & - & - & 0.10 & 8.00  & 268.50  & 0.05 & 1621.72 &
       $0.100\pm0.058$ & $0.414\pm0.445$ & $4.185\pm2.695$ \\
       & - & - & 0.15 & 10.13 & 389.52  & 0.05 & 1621.72 & $0.087\pm0.055$ & $0.175\pm0.147$ & $1.626\pm0.931$\\
       & - & - & 0.25 & 15.24 & 639.69  & 0.05 & 1621.72 & $0.053\pm0.004$ & $0.079\pm0.008$ & $0.682\pm0.479$\\
       & - & - & 0.30 & 18.16 & 758.09  & 0.05 & 1621.72 & $0.051\pm0.004$ & $0.075\pm0.007$ & $0.402\pm0.167$\\
       & - & - & 0.50 & 28.77 & 1257.95 & 0.05 & 1621.72 & $0.050\pm0.003$ & $0.070\pm0.008$ & $0.257\pm0.074$\\
       & - & - & 1.00 & 54.39 & 2483.21 & 0.05 & 1621.72 &  $0.047\pm0.003$ & $0.071\pm0.008$ & $0.140\pm0.034$\\
\bottomrule
\end{tabular}
\end{adjustbox}
\caption{Numerical results comparison of the models for 1D tasks.}
\label{table:efficiency_1d}
\end{table}

\begin{table}[h!]
\begin{adjustbox}{width={\textwidth},totalheight={\textheight},keepaspectratio}
\centering
\begin{tabular}{cccccccccc}
\toprule 
         & & & & \multicolumn{2}{c}{Training}    & \multicolumn{2}{c}{Inference} & \multicolumn{2}{c}{$\epsilon_u$($\times 10^{-2}$)} \\
    \cmidrule(lr){5-6}\cmidrule(lr){7-8}\cmidrule(lr){9-10}         
           & $k$ & $m$ & $p$ & Time(ms)  & Memory(Mb) & Time(ms)  & Memory(Mb) & Poisson 2D & Darcy 2D \\
\midrule
GreenMGNet & 1 & 0 & 0.0664 & 245.78 & 11159.77 & 0.48 & 7354.09 & $1.157\pm0.095$ & $1.433\pm0.079$ \\
           & 1 & 1 & 0.0696 & 270.54 & 11165.76 & 3.23 & 7909.02 &  $0.958\pm0.088$ & $1.272\pm0.089$ \\
           & 1 & 3 & 0.0948 & 377.50 & 11546.68 & 3.78 & 11055.36 &  $0.821\pm0.115$ & $1.175\pm0.098$ \\
           & 1 & 5 & 0.1410 & 559.25 & 17678.45 & 5.34 & 16504.09 & $0.780\pm0.120$ & $1.156\pm0.120$ \\
\midrule
           & 2 & 0 & 0.0047 & 29.50 & 805.61 & 0.64 & 528.51 & $2.956\pm0.063$ & $3.123\pm0.038$ \\
           & 2 & 1 & 0.0086 & 74.50 & 1350.37 & 4.97 & 1263.87 & $1.326\pm0.091$ & $1.584\pm0.086$ \\
           & 2 & 3 & 0.0397 & 194.97 & 5616.91 & 5.77 & 5131.79 & $0.860\pm0.120$ & $1.237\pm0.146$ \\
           & 2 & 5 & 0.0955 & 417.21 & 12845.02 & 7.55 & 11676.11 & $0.820\pm0.178$ & $1.165\pm0.123$ \\
\midrule
           & 3 & 0 & 0.0004 & 27.08 & 80.25 & 0.88 & 51.14 &  $12.98\pm0.011$ & $13.29\pm0.010$ \\
           & 3 & 1 & 0.0045 & 80.47 & 910.93 & 6.81 & 824.52 & $1.964\pm0.041$ & $2.168\pm0.025$ \\
           & 3 & 3 & 0.0369 & 205.38 & 5324.63 & 7.70 & 4839.63 & $0.924\pm0.073$ & $1.241\pm0.057$ \\
           & 3 & 5 & 0.0942 & 427.09 & 12722.09 & 9.49 & 11553.21 & $0.914\pm0.099$ & $1.210\pm0.057$ \\
\midrule
GL & - & - & 0.0010 & 11.84 & 993.78 & 0.21 & 27304.10 & $64.15\pm14.38$ & $63.08\pm13.17$ \\
   & - & - & 0.0050 & 21.35 & 1668.02 & 0.21 & 27304.10 &  $17.00\pm4.419$ & $17.45\pm3.155$ \\
   & - & - & 0.0300 & 102.82 & 5825.12 & 0.21 & 27304.10 & $2.158\pm0.604$ & $2.804\pm0.694$ \\
   & - & - & 0.0500 & 169.74 & 9182.74 & 0.21 & 27304.10 &$1.401\pm0.268$ & $1.837\pm0.309$ \\
   & - & - & 0.1000 & 333.26 & 17531.88 & 0.21 & 27304.10 &$1.037\pm0.129$ & $1.407\pm0.072$ \\
   & - & - & 0.1500 & 498.63 & 25879.83 & 0.21 & 27304.10 &$0.902\pm0.096$ & $1.256\pm0.108$ \\
   & - & - & 0.2000 & 668.88 & 34270.73 & 0.21 & 27304.10 &$0.867\pm0.170$ & $1.202\pm0.107$ \\
\midrule
GL-aug & - & - & 0.0010 & 8.60 & 993.78 & 0.21 & 27362.86 & $60.47\pm19.77$ & $61.84\pm18.84$ \\
 & - & - & 0.0050 & 18.85 & 1668.02 & 0.21 & 27362.86 & $15.15\pm2.491$ & $17.60\pm5.302$ \\
 & - & - & 0.0300 & 102.31 & 5825.12 & 0.21 & 27362.86 &
 $2.639\pm0.670$ & $3.084\pm0.456$ \\
 & - & - & 0.0500 & 169.83 & 9182.74 & 0.21 & 27362.86 & $1.509\pm0.109$ & $2.005\pm0.369$\\
 & - & - & 0.1000 & 336.58 & 17531.88 & 0.21 & 27362.86 &
 $0.928\pm0.061$ & $1.311\pm0.051$\\
 & - & - & 0.1500 & 504.16 & 25879.83 & 0.21 & 27362.86 &
 $0.883\pm0.061$ & $1.238\pm0.091$\\
 & - & - & 0.2000 & 676.43 & 34270.73 & 0.21 & 27362.86 & 
 $0.828\pm0.096$ & $1.185\pm0.089$ \\
\bottomrule
\end{tabular}
\end{adjustbox}
\caption{Numerical results comparison of the models for 2D tasks.}
\label{table:efficiency_2d}
\end{table}

\end{document}